\documentclass[amsppt,11pt]{amsart}
\usepackage{curves}
\usepackage{amsmath}

\newtheorem {theorem}{Theorem}[section]

\newtheorem {lemma}[theorem]{Lemma}
\newtheorem {corollary}[theorem]{Corollary}
\newcounter{conjecture}
\newcounter{remark}

\newcommand{\eqnsection}{
   \renewcommand{\theequation}{\thesection.\arabic{equation}}
   \makeatletter
   \csname @addtoreset\endcsname{equation}{section}
   \makeatother}

\newcommand{\be}{{\begin{equation}}}
\newcommand{\ee}{{\end{equation}}}
\def \bt{\begin{theorem}}
\def \et{\end{theorem}}
\def \bea{\begin{eqnarray}}
\def \eea{\end{eqnarray}}
\def \bas{\begin{eqnarray*}}
\def \eas{\end{eqnarray*}}
\def \bbmod{ {\overline{\bmod} \,}}

\def \al{\alpha}
\def \bb{\beta}
\def \ga{\gamma}

\def \de{\delta}

\def \ep{\epsilon}

\newcommand{\eps}{\varepsilon}

\def \la{\lambda}

\def \om{\omega}

\def \si{\sigma}

\def \th{\theta}

\def \ze{\zeta}

\def \ff{\infty}

\def \wt{\widetilde}

\def \rar{\rightarrow}

\newcommand{\ls}[1]
   {\dimen0=\fontdimen6\the\font \lineskip=#1\dimen0
\advance\lineskip.5\fontdimen5\the\font \advance\lineskip-\dimen0
\lineskiplimit=.9\lineskip \baselineskip=\lineskip
\advance\baselineskip\dimen0 \normallineskip\lineskip
\normallineskiplimit\lineskiplimit \normalbaselineskip\baselineskip
\ignorespaces }
\newcommand{\req}[1]{(\ref{#1})}

\def \R{{\Bbb{R}}}

\def \Z{{\Bbb{Z}}}

\def \AA{{\mathcal A}}

\def \CC{{\mathcal C}}

\def \EE{{\mathcal E}}

\def \GG{{\mathcal G}}
\def \HH{{\mathcal H}}

\def \LL{{\mathcal L}}
\def \MM{{\mathcal M}}
\def \NN{{\mathcal N}}

\def \PP{{\mathcal P}}

\def \RR{{\mathcal R}}
\def \sss{{\mathcal S}}
\def \TT{{\mathcal T}}

\def \VV{{\mathcal V}}

\def \({\left(}
\def \){\right)}

\def \lc{\left\{}
\def \rc{\right\}}

\def \nn{\nonumber}

\def \bc{\begin{center} }
\def \ec{\end{center} }
\def\Bbb{\mathbb}

\begin{document}

\eqnsection
\newcommand{\Ini}{{I_{n,i}}}
\newcommand{\reals}{{\Bbb{R}}}
\newcommand{\F}{{\mathcal F}}
\newcommand{\D}{{\mathcal D}}
\newcommand{\Fn}{{{\mathcal F}_n}}
\newcommand{\Gn}{{{\mathcal G}_n}}
\newcommand{\Hn}{{{\mathcal H}_n}}
\newcommand{\Fp}{{{\mathcal F}^p}}
\newcommand{\Gp}{{{\mathcal G}^p}}
\newcommand{\PPP}{{\mathbf P}}
\newcommand{\Pop}{{P\otimes \PPP}}
\newcommand{\hm}{\HH^\varphi}
\newcommand{\nuw}{{\nu^W}}
\newcommand{\ths}{{\theta^*}}
\newcommand{\beq}[1]{\begin{equation}\label{#1}}
\newcommand{\eeq}{\end{equation}}
\newcommand{\integers}{{\rm I\!N}}
\newcommand{\DT}{{D_{\Bbb T^2}}}
\newcommand{\pDT}{{\partial \DT}}
\newcommand{\E}{{\Bbb E}}
\newcommand{\te}{{\tilde{\delta}}}
\newcommand{\tI}{{\tilde{I}}}
\newcommand{\epn}{{\ep_n}}
\def\var{{\rm Var}}
\def\cov{{\rm Cov}}
\def\one{{\bf 1}}
\def\leb{{\mathcal L}eb}
\def\Ho{{\mbox{\sf H\"older}}}  
\def\thi{{\mbox{\sf Thick}}}
\def\cthi{{\mbox{\sf CThick}}}
\def\late{{\mbox{\sf Late}}}
\def\clate{{\mbox{\sf CLate}}}
\def\plate{{\mbox{\sf PLate}}}
\newcommand{\ffrac}[2]
  {\left( \frac{#1}{#2} \right)}
\newcommand{\calF}{{\mathcal F}}
\newcommand{\dfn}{\stackrel{\triangle}{=}}
\newcommand{\beqn}[1]{\begin{eqnarray}\label{#1}}
\newcommand{\eeqn}{\end{eqnarray}}
\newcommand{\oo}{\overline}
\newcommand{\uu}{\underline}
\newcommand{\bfcdot}{{\mbox{\boldmath$\cdot$}}}
\newcommand{\Var}{{\rm \,Var\,}}
\def\squarebox#1{\hbox to #1{\hfill\vbox to #1{\vfill}}}
\renewcommand{\qed}{\hspace*{\fill}
            \vbox{\hrule\hbox{\vrule\squarebox{.667em}\vrule}\hrule}\smallskip}
\newcommand{\half}{\frac{1}{2}\:}
\newcommand{\beaa}{\begin{eqnarray*}}
\newcommand{\eeaa}{\end{eqnarray*}}
\newcommand{\calK}{{\mathcal K}}
\def\dimm{{\overline{{\rm dim}}_{_{\rm M}}}}
\def\dimp{\dim_{_{\rm P}}}
\def\htaum{{\hat\tau}_m}
\def\htaumk{{\hat\tau}_{m,k}}
\def\htaumkj{{\hat\tau}_{m,k,j}}

\bibliographystyle{amsplain}

\title[Cover Times for planar Brownian Motion and Random Walks]
{
Cover Times for Brownian Motion and Random Walks in two dimensions}

\author[Amir Dembo\,\, Yuval Peres\,\, Jay Rosen\,\, Ofer Zeitouni]
{Amir Dembo$^*$\,\, Yuval Peres$^\dagger$\,\, Jay Rosen$^\ddagger$\,\,
Ofer Zeitouni$^\S$}

\date{submitted: July, 26, 2001; revised: May, 21, 2003.
\newline\indent
$^*$Research  partially supported by NSF grant \#DMS-0072331.
\newline\indent
$^\dagger$Research partially supported by NSF grant \#DMS-9803597.
\newline\indent
$^\ddagger$Research supported, in part, by grants from the NSF  and from
PSC-CUNY.
\newline\indent
$^\S$
The research of all authors was supported, in part, by a US-Israel BSF grant.
}

\begin{abstract}
\noindent
Let $\TT(x,\eps)$ denote the first hitting time of the disc of radius
$\eps$ centered at $x$ for Brownian
motion on the two dimensional torus $\Bbb{T}^2$.
We prove that $\sup_{x\in \Bbb{T}^2}
\TT(x,\eps)/|\log \eps|^2 \to 2/\pi$
as $\eps \rar 0$.  The same applies to Brownian motion on any smooth,
compact connected,  two-dimensional,
Riemannian manifold with unit area and no boundary.
As a consequence, we prove a conjecture, due to Aldous (1989), that
the number of steps it takes a simple random walk to cover all points of
the lattice torus $\Z_n^2$ is asymptotic to $4n^2(\log n)^2/\pi$.
Determining these asymptotics is an essential step toward
analyzing the fractal structure of the set of uncovered sites
before coverage is complete; so far, this structure was 
only studied non-rigorously in the physics literature.
We also establish a conjecture, due to
Kesten and R\'{e}v\'{e}sz, that describes
the asymptotics for the number of steps needed by simple random walk
in $\Z^2$ to cover the disc of radius $n$.
\end{abstract}

\maketitle
\section{Introduction}
In this paper, we introduce a unified method for analyzing
cover times for random walks and Brownian motion in two dimensions,
and resolve
several open problems in this area.

\subsection{Covering the discrete torus}
The time it takes a random walk to cover a finite graph
is a parameter that has been studied intensively by probabilists,
combinatorialists and computer scientists, due to its intrinsic appeal
and its applications to designing universal
traversal sequences~\cite{Ale, bro, bri},
testing graph connectivity \cite{Ale, kar}, and
protocol testing \cite{mihail}; see \cite{aldous2} for
an introduction to cover times.
Aldous and Fill~\cite[Chapter 7]{Ald-F}
consider the cover time for
 random walk on the discrete $d$-dimensional torus
$\Z_n^d=\Z^d/n\Z^d$, and write:
\begin{quote}\it  ``Perhaps surprisingly, the case $d=2$ turns out to be
  the hardest of all explicit graphs for the purpose of estimating
  cover times''. \rm
\end{quote}
The problem of determining the expected cover time $\TT_n$ for $\Z_n^2$
was posed informally
by Wilf \cite{wilf} who called it ``the white screen problem''
and wrote
\begin{quote}\it  ``Any mathematician will want to know how long,
on the average, it takes until each pixel is visited.'' \rm
\end{quote}
(see also \cite[Page 1]{Ald-F}).

In 1989, Aldous \cite{aldous1} conjectured that
$\TT_n/(n\log n)^2\rar 4/\pi$. Aldous noted that the upper bound
$\TT_n/(n\log n)^2 \le 4/\pi+\small{o}(1)$ was
easy, and pointed out the difficulty of obtaining a corresponding lower bound.
A lower bound of the correct order of magnitude was obtained by
Zuckerman~\cite{zuck}, and in 1991,
Aldous~\cite{aldous3} showed that
$\TT_n/\E(\TT_n) \to 1$ in probability.
The best lower bound prior to the present work
is due to Lawler~\cite{Lawler}, who showed that
 $\liminf \E(\TT_n)/(n\log n)^2 \ge 2/\pi$.

Our main result in the discrete setting, is the proof of Aldous's conjecture:
\begin{theorem}
\label{theo-2lattice}
If $\TT_n$ denotes the time it takes
for the simple random walk in $\Z_n^2$ to completely cover $\Z_n^2$, then
\begin{equation}
\lim_{n\to \ff}{\TT_n\over (n\log n)^2}=
\frac{4}{\pi}
\hspace{.1in}\mbox{in probability.}
\end{equation}
\end{theorem}
The main interest in this result is not the value of
the constant, but rather that establishing a limit theorem,
with matching upper and lower bounds, forces one to develop
insight into the delicate process of coverage, and to understand
the fractal structure, and spatial correlations,
 of the configuration of uncovered sites  in $\Z_n^2$ before coverage
 is complete.

The fractal structure of the uncovered set in $\Z_n^2$
has attracted the interest of physicists,
(see \cite{nemi}, \cite{brum-hil} and the references therein), who
used simulations and non-rigorous heuristic arguments to study it.
One cannot begin the rigorous study of this fractal structure without knowing
precise asymptotics for the cover time; an estimate of cover time up to a
bounded factor will not do. See \cite{uncov} for quantitative
results on the uncovered set, based on the ideas of the present paper.

Our proof of Theorem \ref{theo-2lattice} is based on
strong approximation of random walks by Brownian
paths, which reduces that theorem to a question about Brownian
motion on the 2-torus.

\subsection{Brownian motion on surfaces}
For $x$ in the two-dimensional torus
$\Bbb{T}^2$, denote by $\DT (x,\ep)$ the disk of radius $\ep$ centered
at $x$, and consider the hitting time
\[\TT(x,\eps)=\inf \{t>0\,|\,X_t\in \DT (x,\ep)\}.\]
Then
\[\CC_\ep=\sup_{x \in \Bbb{T}^2}\TT(x,\eps) \]
is the $\ep$-covering time of the torus $\Bbb{T}^2$, i.e. the amount of
time needed for the
 Brownian motion $X_t$ to come within $\ep$ of each point in $\Bbb{T}^2$.
Equivalently,
$\CC_\ep$ is the amount of time needed for the Wiener sausage of radius
$\ep$ to completely cover $\Bbb{T}^2$.
We can now state the continuous analog of Theorem \ref{theo-2lattice},
which is the key to its proof.
\begin{theorem}
\label{theo-1p}
For Brownian motion in $\Bbb{T}^2$,
\begin{equation}
\lim_{\ep\rar 0}\frac{\CC_\ep}{
\left(\log \eps\right)^2}=\frac{2}{\pi}
\hspace{.6in}\mbox{a.s.}
\label{p.10}
\end{equation}
\end{theorem}

Matthews~\cite{matt} studied the $\epsilon$-cover time for Brownian
motion on a $d$ dimensional sphere (embedded in $\R^{d+1}$)
and on a $d$-dimensional projective space
(that can be viewed as the quotient of the sphere
 by reflection). He calls these questions the ``one-cap problem''
 and ``two-cap problem'', respectively. Part of the motivation for
 this study is a technique for viewing multidimensional data
developed by Asimov~\cite{Asimov}. Matthews obtained sharp
asymptotics for all dimensions $d \ge 3$, but for the
more delicate two dimensional
case, his upper and lower bounds had a ratio of $4$ between them;
he conjectured the upper bound was sharp.
We can now resolve this conjecture; rather than handling each surface
separately, we establish the following extension of Theorem
\ref{theo-1p}. See Section \ref{sec-manifold} for definitions
and references concerning Brownian motion on manifolds.

\begin{theorem}
\label{theo-m1p}
Let $M$ be a smooth, compact, connected two-dimensional,
Riemannian manifold without boundary.
Denote by $\CC_\ep$ the $\ep$-covering time of $M$, i.e.,
 the amount of time needed for the
 Brownian motion to come within
(Riemannian) distance $\ep$ of each point in $M$.
Then
\begin{equation}
\lim_{\ep\rar 0}\frac{\CC_\ep}{
\left(\log \eps\right)^2}=\frac{2}{\pi} A \quad a.s.,
\label{mp.10}
\end{equation}
where $A$ denotes the Riemannian area of $M$.
\end{theorem}
(When $M$ is a sphere,
this indeed corresponds to the upper bound in \cite{matt},
 once a computational error in \cite{matt} is corrected;
  the hitting time in (4.3) there is twice what it should be.
This error led to doubling the upper and the lower bounds for cover
 time in \cite[Theorem 5.7]{matt}).

\subsection{Covering a large disk by random walk in $\Z^2$}

Over ten years ago, Kesten (as quoted by Aldous~\cite{aldous1}
and Lawler~\cite{Lawler})
 and R\'{e}v\'{e}sz \cite{Rev} independently considered a problem
about simple random walks in $\Z^2$:
\em How long does it take for the walk to
completely cover
the disc of radius $n$?
\rm Denote this time by $T_n$.
Kesten and R\'{e}v\'{e}sz proved that
\begin{equation}
\label{KR1}
e^{-b/t} \le
\liminf_{n\to\infty} \PPP(\log T_n\leq t(\log n)^2) \le
\limsup_{n\to\infty} \PPP(\log T_n\leq t(\log n)^2) \le
e^{-a/t}.
\end{equation}
for certain  $0<a<b<\ff$.
R\'{e}v\'{e}sz~\cite{Rev}
conjectured that the limit exists and has the form $e^{-\lambda/t}$
for some (unspecified) $\lambda$.
Lawler~\cite{Lawler}  obtained (\ref{KR1})
with the constants $a=2,\,b=4$  and quoted a conjecture of Kesten
that the limit equals $e^{-4/t}$. We can now prove this:
\begin{theorem}
\label{theo-4a}
If $T_n$ denotes the time it takes
for the simple random walk in $\Z^2$ to completely cover the disc of radius
$n$, then
\begin{equation}
\label{eq-KRa}
\lim_{n\to\infty} \PPP(\log T_n\leq t(\log n)^2) =e^{-4/t}.
\end{equation}
\end{theorem}

\subsection{A birds-eye view}
The basic approach of this paper, as in \cite{DPRZ4}, is to
control $\ep$-hitting times using
excursions between concentric circles.
The number of excursions between two fixed concentric circles
before $\epsilon$-coverage
is so large, that the $\ep$-hitting times will
necessarily be concentrated near their conditional means
given the excursion counts (see Lemma~\ref{lem1}).

The key idea in the proof of the lower bound in Theorem
\ref{theo-1p},
is to  control excursions on many scales simultaneously,
leading to a `multi-scale refinement' of the classical second
moment method.
This is inspired by techniques from probability on trees,
in particular the analysis of first-passage percolation by Lyons and
Pemantle~\cite{LP}. The approximate tree structure that we
(implicitly) use arises by considering circles  of varying radii
around different centers;
 for fixed centers $x,y$, and ``most'' radii $r$ (on a logarithmic
scale) the discs $\DT (x,r)$ and $\DT (y,r)$ are either well-separated
(if $r \ll d(x,y)$) or almost coincide (if $r \gg d(x,y)$).
 This tree structure was also the key to our work in \cite{DPRZ4},
but the dependence problems encountered in the present work are more
severe.  While in \cite{DPRZ4} the number of
macroscopic excursions was bounded, here it is large;
In the language of trees, one
can say that while in \cite{DPRZ4} we studied the maximal number of visits
to a leaf until visiting the root, here we study the number of visits to
the root until every leaf has been visited. For the analogies between
trees and Brownian excursions to be valid, the effect of the initial
and terminal points of individual excursions must be controlled.
To prevent conditioning on the endpoints of the numerous
macroscopic excursions to
affect the estimates, the ratios between radii of  even the largest pair
of concentric circles where excursions are counted, must grow to infinity
as $\epsilon$ decreases to zero.

Section \ref{sec-hit} provides
simple lemmas
which will be useful in exploiting
the link between  excursions and $\ep$-hitting times.
These lemmas are then used
to obtain the upper bound in Theorem
\ref{theo-1p}. In Section \ref{sec-KRlowerbound}
we explain how to obtain the analogous lower bound,
leaving some technical details to
lemmas which are proven in Sections \ref{1momest}-\ref{2momest}.
In Section
\ref{sec-etc} we prove the lattice torus covering time conjecture,
Theorem \ref{theo-2lattice}, and in Section \ref{RKproof} we prove the
Kesten-R\'{e}v\'{e}sz conjecture, Theorem \ref{theo-4a}.
In Section \ref{sec-manifold} we consider
Brownian motion on manifolds and prove Theorem
\ref{theo-m1p}.
Complements and open problems are collected in the final section.

\section{Hitting time estimates and upper bounds}\label{sec-hit}
\label{sec-upperbound}
We start with some definitions.
Let $\{W_t\}_{t \ge 0}$ denote  planar Brownian motion started at
the origin. We use
$\Bbb{T}^2$
to denote the two dimensional torus, which we identify
with  the set $(-1/2,1/2]^2$.
The distance between $x,y \in \Bbb{T}^2$, in the natural metric,
is denoted $d(x,y)$.
Let $X_t=W_t \, \bbmod \,\Z^2$
denote the Brownian motion on $\Bbb{T}^2$,
where
$a\,\bbmod \,\Z^2=[a +(1/2,1/2)] \bmod \,\Z^2 - (1/2,1/2)$.
Throughout,
$D(x,r)$
and $\DT (x,r)$ denote the open discs
of radius $r$ centered at $x$,
in $\reals^2$ and in $\Bbb{T}^2$, respectively.

Fixing $x\in \Bbb{T}^2$ let
$\tau_\xi=\inf\{t \geq 0\, :\, X_t \in \pDT (x,\xi) \}$ for $\xi>0$.
Also let $\wt{\tau}_\xi = \inf \{t \geq 0 : B_t \in \partial
D(0,\xi) \}$,
for a standard Brownian motion $B_t$ on $\reals^2$.
For any $x \in \Bbb{T}^2$, the natural bijection
$i=i_x:\DT (x,1/2) \mapsto D(0,1/2)$
with $i_x(x)=0$ is an isometry, and for any
$z \in \DT (x,1/2)$ and
Brownian motion $X_t$ on $\Bbb{T}^2$ with $X_0=z$,
we can find a Brownian motion
$B_t$ starting
at $i_x(z)$ such that $\tau_{1/2}=\wt{\tau}_{1/2}$ and
$\{ i_x(X_t), t \leq \tau_{1/2} \} = \{ B_t,
t \leq \wt{\tau}_{1/2} \}$.
We shall hereafter use $i$ to denote $i_x$, whenever
the precise value of $x$ is understood from the context, or 
does not matter.

We start with some uniform estimates on the hitting times
$\E^y(\tau_r)$.
\begin{lemma}
\label{lem-hit}
For some $c<\infty$ and all $r>0$ small enough,
\beq{tohelet-max}
\|\tau_r\|:=
\sup_{y}  \E^y (\tau_r) \leq c |\log r|  \,.
\end{equation}
Further, there exists $\eta(R) \to 0$ as $R \to 0$, such that
for all $0<2r \leq R$, $x \in \Bbb{T}^2$,
\beqn{tohelet0}
{(1 -\eta)  \over \pi}\log \ffrac{R}{r} &\leq&
\inf_{y \in \pDT (x,R)}  \E^y (\tau_r) \nonumber
\\
&\leq& \sup_{y \in \pDT (x,R)}  \E^y (\tau_r) \leq
{(1 +\eta)  \over \pi}\log \ffrac{R}{r} \,.
\eeqn
\end{lemma}

\noindent
{\bf Proof of Lemma~\ref{lem-hit}:}
Let $\Delta$ denote the Laplacian, which on $\Bbb{T}^2$
is just the Euclidean Laplacian with periodic boundary conditions.
It is well known that
for any $x\in \Bbb{T}^2$ there exists
a
Green's function
$G_x(y)$, defined for $y \in \Bbb{T}^2 \setminus \{x\}$,
such that $\Delta G_x = 1$ and
$F(x,y)=G_x(y) + \frac{1}{2\pi} \log d(x,y)$ is continuous
on $\Bbb{T}^2 \times \Bbb{T}^2$
(c.f. \cite[p. 106]{Aubin} or \cite{Eells-Sampson} where this is
shown in the more general context of smooth, compact
two-dimensional Riemannian manifold without boundary).
For completeness, we explicitly construct such $G_x(\cdot)$ at
the end of the proof.

Let $e(y)=E^y(\tau_r)$.
We have Poisson's equation ${1\over 2}\Delta e=-1$
on $\Bbb{T}^2 \setminus \DT (x,r)$ and $e=0$ on $\pDT (x,r)$.
Hence, with $x$ fixed,
\begin{equation}
\Delta \(G_x+{1\over 2}e\)=0\hspace{.2in}\mbox{on}
\hspace{.1in} \Bbb{T}^2 \setminus \DT (x,r).
\label{k2}
\end{equation}
Applying the maximum principle for the harmonic function
$G_x+{1\over 2}e$ on $\Bbb{T}^2 \setminus \DT (x,r)$, we see
that for all $y\in \Bbb{T}^2 \setminus \DT (x,r)$,
\begin{equation}
\inf_{z\in \pDT (x,r)} G_x(z) \leq G_x(y)
+{1\over 2}e(y)\leq \sup_{z\in \pDT (x,r) }G_x(z).
\label{k3}
\end{equation}
Our lemma follows then, with
\beaa
\eta(R) &=& \frac{2\pi}{\log 2} \sup_{x \in \Bbb{T}^2} \;
                                \sup_{y,z\in \DT (x,R)}|F(x,z)-F(x,y)|\\
c & = & (1/\pi) + [(1/\pi) \log {\rm diam(\Bbb{T}^2)}
+4\sup_{x,y\in \Bbb{T}^2} |F(x,y)|]/\log 4 < \infty \;,
\eeaa
except that we have proved \req{tohelet-max} so far only for
$y \notin \DT (x,r)$. To complete the proof, fix $x' \in \Bbb{T}^2$ with
$d(x,x')=3 \rho>0$.
For $r<\rho$, starting at $X_0=y \in \DT (x,r)$,
the process $X_t$ hits $\pDT (x,r)$ before it hits $\pDT (x',r)$.
Consequently, $E^y(\tau_r) \leq c |\log r|$
also for such $y$ and $r$,
establishing \req{tohelet-max}.

Turning to construct $G_x(y)$, we use the representation
$\Bbb{T}^2=(-1/2,1/2]^2$. Let $\phi \in C^\ff(\reals)$ be such
that $\phi=1$ in a small neighborhood of $0$, and $\phi=0$
outside a slightly larger neighborhood of $0$. With $r=|z|$ for
$z=(z_1,z_2)$, let
\[
h(z)=-{1\over 2\pi}\phi(r)\log r
\]
and note that by Green's theorem
\begin{equation}
\int_{\Bbb{T}^2} \Delta h(z)\,dz=1.\label{gf.1}
\end{equation}
Recall that for any function $f$
which depends only on $r=|z|$
\[\Delta f=f''+{1\over r}f',\]
 and therefore, for $r>0$
$$
\Delta h(z)=-{1\over 2\pi}(\phi''(r)\log r+{2+\log r\over r}\phi'(r)).
$$
Because of the support properties of $\phi(r)$ we see that
$H(z)=\Delta h(z) - 1$ is a $C^\ff$
function on $\Bbb{T}^2$, and consequently has an expansion in Fourier series
$$
H(z)=\sum_{j,k=0}^\ff a_{j,k}\cos (2\pi j z_1)\cos (2\pi k z_2)
$$
with $a_{j,k}$ rapidly decreasing. Note that as a consequence of
(\ref{gf.1}) we have $a_{0,0}=0$. Set
$$
F(z)=\sum_{\stackrel{j,k=0}{(j,k)\neq (0,0)}}^\ff {a_{j,k}\over
4\pi^2(j^2+k^2)}
\cos (2\pi j z_1)\cos (2\pi k z_2).
$$
The function $F(z)$ is then a $C^\ff$ function on $\Bbb{T}^2$ and
it satisfies $\Delta F=-H$. Hence, if we set
$g(z)=h(z)+F(z)$ we have $\Delta g(z)=1$ for $|z|>0$
and $g(z)+{1\over 2\pi}\log |z|$ has a continuous
extension to all of $\Bbb{T}^2$.
The
Green's function for $\Bbb{T}^2$  is then
$G_x(y)=g( (x-y)_{\Bbb{T}^2} )$.
\qed

Fixing $x\in \Bbb{T}^2$ and constants $0<2 r \leq R<1/2$ let
\begin{equation}
\tau^{(0)}=\inf \{t\geq 0\,|\,X_t\in \pDT (x,R) \}\label{c2.1}
\end{equation}
\begin{equation}
\si^{(1)}=\inf\{t\geq 0\,|\,X_{t+\tau^{(0)}}\in
\pDT (x,r) \}\label{c2.2}
\end{equation}
and define inductively for $j=1,2,\ldots$
\begin{equation}
\tau^{(j)}=\inf \{t\geq \si^{(j)}\,|\,X_{t+\frak{T}_{j-1}}
\in \pDT (x,R)  \},\label{c2.3}
\end{equation}
\begin{equation}
\si^{(j+1)}=\inf\{t\geq 0\,|\,X_{t+\frak{T}_j}\in \pDT (x,r) \},
\label{c2.4}
\end{equation}
where $\frak{T}_{j}=\sum_{i=0}^{j}\tau^{(i)}$ for $j=0,1,2,\ldots$.
Thus, $\tau^{(j)}$ is the length of the $j$-th excursion $\EE_j$ from
$\pDT (x,R) $ to itself via $\pDT (x,r) $,
and $\si^{(j)}$ is the amount of time it takes to hit $\pDT (x,r)$
during the $j$-th excursion $\EE_j$.

The next lemma, which shows that excursion times are concentrated around
their mean, will be used to relate excursions to hitting times.

\begin{lemma}
\label{lem-ld}
With the above notation, for any
$N \geq N_0$, $\de_0>0$ small enough, $0<\de<\de_0$,
$0<2 r \leq R<R_1(\delta)$, and $x,x_0 \in \Bbb{T}^2$,
\begin{equation}
\PPP^{x_0}\(\sum_{j=0}^N \tau^{(j)}\leq (1-\de)N {1 \over \pi} \log(R/r)
\)\leq e^{-C\de^2 N}\label{c2.10a}
\end{equation}
and
\begin{equation}
\PPP^{x_0}\(\sum_{j=0}^N \tau^{(j)}\geq (1+\de) N {1 \over \pi} \log(R/r)
\)\leq e^{-C\de^2 N}\label{c2.10b}
\end{equation}
Moreover, $C=C(R,r)>0$ depends only upon
$\delta_0$ as soon as $R>r^{1-\delta_0}$.
\end{lemma}

\noindent
{\bf Proof of Lemma~\ref{lem-ld}:}
Applying Kac's moment formula for the first
hitting time $\tau_r$ of the strong Markov process $X_t$
(see \cite[Equation (6)]{Fitzsimmons-Pitman}), we see that
for any $\th<1/\|\tau_r\|$,
\begin{equation}
\sup_{y}  \E^y (e^{\th\tau_r})\leq {1\over 1-\th\|\tau_r\|} \,.
\label{k6}
\end{equation}
Consequently, by (\ref{lem-hit}) we have that for some $\la>0$,
\begin{equation}
\label{expmomaj}
\sup_{0< r \leq r_0} \sup_{x,y} \E^{y}(e^{\la \tau_r/|\log r|})<\ff.
\end{equation}
By the strong Markov property of $X_t$ at $\tau^{(0)}$
and at $\tau^{(0)}+\si^{(1)}$ we then deduce that
\begin{equation}
\label{expmoma}
\sup_{0<2 r \leq R< r_0} \sup_{x,y} \E^{y}(e^{\la \frak{T}_1/|\log r|})<\ff.
\end{equation}

Fixing $x\in \Bbb{T}^2$ and $0<2 r \leq R<1/2$ let
$\tau=\tau^{(1)}$ and $v=\frac{1}{\pi}\log(R/r)$.
Recall that $\{X_t: t \leq \tau_R \}$
starting at $X_0=z$ for some $z \in \pDT (x,r) $,
has the same law as $\{B_t : t \leq \wt{\tau}_R \}$
starting at
$B_0 = i(z) \in \partial D(0,r)$. Consequently,
\beq{amir-may1}
\| \tau_R \|_R
:= \sup_{x} \sup_{z \in \DT (x,R)} \E^z(\tau_R)
\leq
\E^0(\wt{\tau}_R)=\frac{R^2}{2}
 \to_{R \to 0} 0 \,,
\end{equation}
by the
radial symmetry of the Brownian motion
$B_t$.

By the strong Markov property of $X_t$ at $\tau^{(0)}+\si^{(1)}$
we thus have that
$$
\E^y(\tau_r) \leq \E^y(\tau) \leq \E^y(\tau_r)+\|\tau_R\|_R
\quad
\quad \forall y \in \pDT (x,R)
$$
Consequently, with $\eta=\delta/6$, let
$R_1(\delta) \leq r_0$
be small enough
so that \req{tohelet0} and \req{amir-may1} imply
\beqn{tohelet}
(1 -\eta) v &\leq&
\inf_{x} \inf_{y \in \pDT (x,R) }  \E^y (\tau) \nonumber \\
&\leq&
 \sup_{x} \sup_{y \in \pDT (x,R) }  \E^y (\tau) \leq
(1 +2\eta) v \,,
\eeqn
whenever $R \leq R_1$.
It follows from
\req{expmoma}
and
\req{tohelet}
that there exists a universal
constant $c_4<\infty$ such that for $\rho=c_4 |\log r|^2$
and all $\theta \geq 0$,
\beqn{new-star1}
&& \sup_{x} \sup_{y \in \pDT (x,R) } \E^y (e^{-\theta \tau})
\nn \\
&\leq&
 1 - \theta \inf_{x}\inf_{y \in \pDT (x,R)} \E^y (\tau) + \frac{\theta^2}{2}
\sup_{x}\sup_{y \in \pDT (x,R) } \E^y (\tau^2)
\nn \\
&\leq&
 1 - \theta (1-\eta) v + \rho \theta^2
\leq \exp(\rho \theta^2 - \theta (1-\eta) v)
\eeqn
Since $\tau^{(0)} \geq 0$, using Chebyshev's
inequality we bound the left hand side of \req{c2.10a} by
\beqn{new-star2}
\PPP^{x_0}\Bigl(\sum_{j=1}^N \tau^{(j)}\leq (1-6\eta) v N\Bigr)
&\leq&
e^{\theta (1-3\eta) v N} \E^{x_0}
\Bigl( e^{-\theta \sum_{j=1}^N \tau^{(j)}} \Bigr) \nn \\
&\leq& e^{-\theta v N \de/3}
\Bigl[ e^{\theta (1-\eta) v}
\sup_{y \in \pDT (x,R)} \E^y (e^{-\theta \tau}) \Bigr]^N  ,
\eeqn
where the last inequality follows by the strong Markov property of $X_t$ at
$\{\frak{T}_j\}$. Combining \req{new-star1} and \req{new-star2} for
$\theta=\delta v / (6 \rho)$, results
in
\req{c2.10a}, where
$C=v^2/36\rho >0$ is bounded below by
$\delta_0^2/(36 c_4 \pi^2)$ if $r^{1-\delta_0}<R$.

To prove (\ref{c2.10b}) we first note that for $\theta=\la/|\log r|>0$
and $\la>0$ as in \req{expmoma}, it follows that
\[
\PPP^{x_0}\( \tau^{(0)}\geq
\frac{\de}{3} v N\)
\leq e^{-\theta v (\de/3) N}\E^{x_0}(e^{\la\tau^{(0)}/|\log r|})
\leq c_5 e^{-c_6 \delta N} \,,
\]
where $c_5<\infty$ is a universal constant and $c_6=c_6(r,R)>0$
does not depend upon $N$, $\delta$ or $x_0$ and is bounded below by
some $c_7(\de_0)>0$ when $r^{1-\de_0} < R$.
Thus, the proof of  (\ref{c2.10b}), in analogy to that of (\ref{c2.10a}),
comes down
to bounding
\begin{equation}
\PPP^{x_0}\Bigl(\sum_{j=1}^N \tau^{(j)} \geq (1+4\eta) v N
\Bigr)
\leq e^{-\th \de vN/3}
\Bigl( e^{-\th (1+2 \eta) v}
\sup_{y \in \pDT (x,R) } \E^y (e^{\theta \tau}) \Bigr)^N
\label{c2.10e}
\end{equation}
Noting that, by
\req{expmoma} and \req{tohelet},
there exists a universal
constant $c_8<\infty$ such that for $\rho=c_8 |\log r|^2$
and all $0<\th < \la/(2 |\log r|)$,
\beaa
\sup_{x} \sup_{y \in \pDT (x,R) } \E^y (e^{\theta \tau})
&\leq& 1 +\th (1+2\eta) v +
\sup_{x} \sup_{y \in \pDT (x,R) }
\sum_{n=2}^\ff \frac{\th^n}{n!} \E^y (\tau^n)
\\
&\leq& 1 +\th (1+2\eta) v + \rho \th^2 \leq
\exp(\th (1+2\eta) v + \rho \th^2) \,,
\eeaa
the proof of (\ref{c2.10b}) now follows as in the
proof of (\ref{c2.10a}).
\qed

\begin{lemma}
\label{lem-hitprob}
For any $\de>0$ we can find $c<\ff$ and  $\eps_0>0$ so that for all
$\eps\leq \eps_0$
 and $y\geq 0$
\begin{equation}
\label{c2.8}
\PPP^{x_0}\(\TT(x,\eps)\geq y (\log \eps)^2\)\leq c\eps^{(1-\de)\pi y}
\end{equation}
for all $x,x_0\in \Bbb{T}^2$.
\end{lemma}

\noindent
{\bf Proof of Lemma~\ref{lem-hitprob}:} We use the notation of the last
lemma and its
proof, with $R<R_1(\de)$ and $r=R/e$ chosen for convenience so that
$\log (R/r)=1$.
Let $n_\eps:= (1-\de)\pi y (\log \eps)^2$.
Then,
\begin{equation}
\PPP^{x_0}\(\TT(x,\eps)\geq y (\log \eps)^2\)\label{c2.9}
\end{equation}
\[
\leq
\PPP^{x_0}\(\TT(x,\eps)\geq \sum_{j=0}^{n_\eps}
\tau^{(j)}\) +
\PPP^{x_0}\(\sum_{j=0}^{n_\eps} \tau^{(j)}\geq y (\log \eps)^2\)\]
It follows from Lemma~\ref{lem-ld} that
\begin{equation}
\PPP^{x_0}\(\sum_{j=0}^{n_\eps} \tau^{(j)}\geq y
(\log \eps)^2\)\leq e^{-C'y(\log \eps)^2}\label{c2.11}
\end{equation}
for some $C'=C'(\de)>0$. On the other hand, the first probability in the
second line of
(\ref{c2.9}) is
bounded above by the probability of
$B_t$ not hitting $i (\DT (x,\ep))=D(0,\ep)$ during
$n_\ep$ excursions,
each starting at $i (\pDT (x,r) ) = \partial D(0,r) $
and ending at $i (\pDT (x,R)) = \partial D(0,R)$,
so that
\begin{equation}
\PPP^{x_0}\(\TT(x,\eps)\geq \sum_{j=0}^{n_\eps} \tau^{(j)}\)\leq
 \(1-{1\over \log \frac{R}{\eps}}\)^{n_\eps}\leq
e^{-(1-\de)\pi y |\log \eps|}\label{c2.12}
\end{equation}
and (\ref{c2.8}) follows.
\qed

We next show that
\begin{equation}
\label{3.1b}
\limsup_{\eps\to0}
\sup_{x\in \Bbb{T}^2}
\frac{\TT(x,\eps)}{\left(\log \eps\right)^2} \le  \frac{2}{\pi}
\,,\quad a.s.
\end{equation}
from which the upper bound for \req{p.10} follows.

Set $h(\ep)=|\log \ep|^2$.
Fix
$\de >0$, and
set $\tilde\ep_n=e^{-n}$ so that
\begin{equation}
h(\tilde\ep_{n+1}) = (1+\frac{1}{n})^{2} h(\tilde\ep_n).
\label{3.2}
\end{equation}
Since, for $\tilde\ep_{n+1}\leq \ep\leq \tilde\ep_n$ we have
\begin{equation}
\label{3.3}
\quad \quad
{\TT(x,\tilde\ep_{n+1})\over
h(\tilde\ep_{n+1})}= {h(\tilde\ep_{n})\over
h(\tilde\ep_{n+1}) } {\TT(x,\tilde\ep_{n+1})\over h(\tilde\ep_{n})} \geq
(1+\frac{1}{n})^{-2}
{\TT(x,\ep)\over h(\ep)} \, .
\end{equation}
Fix $x_0 \in \Bbb{T}^2$ and let $\{ x_j : j=1,\ldots,\bar K_n\}$,
denote a maximal collection of points in
$\Bbb{T}^2$,
such that $\inf_{\ell \neq j}
\, d(x_\ell,x_j)
\geq \delta \tilde\ep_n$.
Let $a=(2+\de)/(1-10\de)$ and
$\AA_{n}$ be the set of $1\leq j\leq \bar K_n$, such that
$$
\TT(x_j,(1-\de)\tilde\ep_n) \geq (1-2\de) a h(\tilde\ep_n)/\pi.
$$
It follows by
Lemma~\ref{lem-hitprob} that
\[
\PPP^{x_0}(\TT(x,(1-\de)\tilde\ep_n) \geq (1-2\de) a
h(\tilde\ep_n)/\pi)
\leq c\, \tilde\ep_n^{\,\,(1-10\delta) a} \,,
\]
for some $c=c(\delta)<\infty$, all sufficiently large $n$ and
any $x \in \Bbb{T}^2$.
Thus, for all sufficiently large $n$, any $j$ and $a>0$,
\begin{equation}
\label{3.3j}
\PPP^{x_0}(j\in \AA_n)\leq c \,\tilde\ep_n^{\,\,(1-10\de) a}\,,
\end{equation}
implying that
$$
\sum_{n=1}^\infty \PPP^{x_0}(|\AA_n|\geq 1)
\leq \sum_{n=1}^\infty \E^{x_0}|\AA_n| \leq c'\sum_{n=1}^\infty
\tilde\ep_n^{\,\,\de}< \infty\,.
$$
By Borel-Cantelli, it follows that $\AA_n$ is empty a.s. for all
$n>n_0(\omega)$ and some $n_0(\om)<\infty$. By (\ref{3.3}) we then have
for some $n_1(\de,\om)<\infty$ and all $n>n_1(\omega)$
$$
\sup_{\ep \leq {\tilde\ep}_{n_1}} \,
\sup_{x\in \Bbb{T}^2}
\frac{\TT(x,\eps)}{
\left( \log \eps \right)^2} \le  \frac{a}{\pi} \,,\,$$
and (\ref{3.1b}) follows by taking $\delta\downarrow 0$.
\qed

\section{Lower bound for covering times}\label{sec-KRlowerbound}
\label{sec-KRprth2}

Fixing $\de>0$ and $a<2$, we prove in this section that
\begin{equation}
\liminf_{\ep\rar 0}\frac{\CC_\ep} {(\log \eps)^2}
\geq (1-\de) \frac{a}{\pi}
\hspace{.6in}\mbox{a.s.}
\label{p.1}
\end{equation}
In view of (\ref{3.1b}), we then obtain Theorem \ref{theo-1p}.

We start by constructing an almost sure lower bound on $\CC_\ep$
for a specific deterministic sequence $\ep_{n,1}$.
To this end, fix $\ep_1 \leq R_1(\de)$ as in Lemma \ref{lem-ld}
and the square $S=[\ep_1,2\ep_1 ]^2$.
Let $\ep_k=\ep_1 (k!)^{-3}$ and $n_k=3a k^2\log k$.
Per fixed $n \geq 3$, let $\ep_{n,k} = \rho_n \ep_n (k!)^3$
for 
$\rho_n=n^{-25}$ 
and $k=1,\ldots,n$.
Observe that $\ep_{n,1}=\rho_n \ep_{n}$, $\ep_{n,n}=\rho_n \ep_1$,
and $\ep_{n,k}\leq \rho_n \ep_{n+1-k} \leq \ep_{n+1-k}$
for all $1\leq k\leq n$.
Recall the natural bijection
$i:\DT (0,1/2) \mapsto D(0,1/2)$.
For any $x \in S$, let
$\RR^x_n$ denote the time
until $X_t$ completes $n_n$ excursions from
$i^{-1} (\partial D(x,\ep_{n,n-1}) )$
to $i^{-1} (\partial D(x,\ep_{n,n}) )$.
(In the notations of Section \ref{sec-hit}, if we set
$R=\ep_{n,n}$ and $r=\ep_{n,n-1}$,
then $\RR^x_n=\sum_{j=0}^{n_n} \tau^{(j)}$).
Note that $i^{-1} (\partial D(x,\ep_{n,k}) )$ is just $\partial \DT (i^{-1}
(x),\ep_{n,k})$, but
the former notation will allow easy generalization to the case of general
manifolds
treated in Section \ref{sec-manifold}.

For $x \in S$, $2\leq k\leq n$ let $N_{n,k}^x$ denote the number of
excursions
of $X_t$
from
$i^{-1} (\partial D(x,\ep_{n,k-1}) )$
to $i^{-1} (\partial D(x,\ep_{n,k}) )$
until time $\RR^x_n$. Thus,
$N^x_{n,n}=n_n=3 a n^2\log n$.
A point $x\in S$ is called {\bf $n$-successful} if
\beq{pmperf}
N^x_{n,2}=0, \qquad n_k-k \leq N^x_{n,k} \leq n_k+k  \quad
\forall k=3,\ldots,n-1\,.
\end{equation}
In particular, if $x$ is $n$-successful, then 
$\TT(i^{-1}(x),\ep_{n,1})>\RR^x_n$.

For $n \geq 3$ we partition $S$ into
$M_n=\ep_1^2/(2 \ep_{n})^2=
(1/4)
\prod_{l=1}^{n }l^6$ non-overlapping
squares of edge length $2 \ep_{n}=2 \ep_1/({n}!)^3$, with
$x_{n,j}$, $j=1,\ldots,M_n$ denoting the centers of these squares.
Let $Y(n,j)$, $j=1,\ldots,M_n$ be the sequence of random variables defined by
\[Y(n,j)=1\hspace{.1in}\mbox{if $x_{n,j}$ is $n$-successful}\]
and $Y(n,j)=0$ otherwise. Set
$\bar{q}_{n}=\PPP(Y(n,j)=1)=\E(Y(n,j))$, noting that this probability is
independent of $j$ (and of the value of $\rho_n$).

The next lemma, which is
a direct consequence of Lemmas \ref{moment1} and \ref{moment2},
provides bounds on the
first and second moments of $Y(n,j)$,
that are used in order to show the existence of
at least one $n$-successful point $x_{n,j}$ for large enough $n$.
\begin{lemma}\label{pmomlb}
There exists $\de_n \to 0$ such that for all $n \geq 1$,
\begin{equation}\qquad
\bar{q}_{n}= \PPP(x\, \mbox{ is $n$-successful}) \geq
\ep_{n}^{a+\de_n} \, \,,
\label{pmomlb.1}
\end{equation}
For some
$C_0<\infty$ and all $n$, if $|x_{n,i}-x_{n,j}| \geq 2 \ep_{n,n}$, then
\begin{equation}
\E(Y(n,i)Y(n,j))\leq (1+C_0 n^{-1} \log n) \bar{q}_n^2 \;.
\label{pm3.50}
\end{equation}
Further, for
any $\ga>0$ we can find $C=C(\ga) <\ff$ so that for all $n$ and
$l=l(i,j) = \max \{k \leq n\,: \,|x_{n,i}-x_{n,j}|
\geq 2 \ep_{n,k} \} \vee 1$,
\begin{equation}
\E(Y(n,i)Y(n,j))\leq \bar{q}_n^2 C^{n-l}
n^{39} \left(\frac{\ep_{n,n}}{\ep_{n,l+1}}\right)^{a+\ga} \;.
\label{pm3.5}
\end{equation}
\end{lemma}

Fix $\ga>0$ such that $2-a-\ga>0$. By (\ref{pmomlb.1})
for all $n$ large enough,
\begin{equation}
\E\(\sum_{j=1}^{M_n} Y(n,j) \)=M_n \bar{q}_{n} \geq \ep_n^{-(2-a-\ga)} .
\label{pm3.7}
\end{equation}
In the sequel, we let $C_m$ denote generic finite constants that are
independent of $n$, $l$, $i$ and $j$.
Recall that there are at most $C_1 \ep_{n,l+1}^2\ep_{n}^{-2}$ points
$x_{n,j}$,
$j \neq i$, in $D(x_{n,i},2\ep_{n,l+1})$.
Further, our choice of $\rho_n$ guarantees that
$(\ep_{n,n}/\ep_n)^2 \leq C_2 M_n 
n^{-50}$.
Hence, it follows from \req{pm3.5} that for $n-1 \geq l \geq 1$,
\bea
V_l &:=& (M_n \bar{q}_n)^{-2}
\sum_{\stackrel{i \neq j=1}{l(i,j)=l}}^{M_n} \E\Big( Y(n,i) Y(n,j) \Big)
\label{jpatch2} 
\\
&\leq& C_1 M_n^{-1} \ep_{n,l+1}^2 \ep_n^{-2}
C^{n-l} n^{39} \ffrac{\ep_{n,l+1}}{\ep_{n,n}}^{-a-\ga}
\nn 
\\
\nn
&\leq& C_1 C_2 n^{-3} C^{n-l} \ffrac{\ep_{n,l+1}}{\ep_{n,n}}^{2-a-\ga}
\;,
\eea
and since $(\ep_{n,l+1}/\ep_{n,n})\leq (\ep_{n-l}/\ep_1)$
for all $1 \leq l \leq n-1$, we deduce that
\beq{p-new}
\sum_{l=1}^{n-1} V_l \leq
C_3  n^{-3} \sum_{j=1}^{\ff} C^j \ep_{j}^{2-a-\ga} \leq C_4 n^{-3}
\,.
\end{equation}
We have, by Chebyshev's inequality (see \cite[Theorem 4.3.1]{Alon})
and \req{pm3.50}, that
\beaa 
\PPP(\sum_{j=1}^{M_n} Y(n,j) = 0)&\leq&
(M_n \bar{q}_n)^{-2}
\E\lc \Big(\sum_{i=1}^{M_n} Y(n,i) \Big)^2  \rc -1  \nn \\
&\leq& (M_n \bar{q}_n)^{-1} + C_0 n^{-1} \log n
 + \sum_{l=1}^{n-1} V_l \,.
\eeaa 
Combining 
this with \req{pm3.7} and \req{p-new}, 
we see that \beq{p-new2}
\PPP(\sum_{j=1}^{M_n} Y(n,j) = 0) \leq C_5 n^{-1} \log n\,.
\end{equation}

The next lemma relates the notion of $n$-successful to
the $\ep_{n,1}$-hitting time.
\begin{lemma}\label{lem1}
For each $n$ let $\VV_n$ be a finite subset of $S$ with cardinality
bounded by $e^{o(n^2)}$.  There exists $m(\om)<\ff$ a.s. such that for all
$n\geq m$ and all $x\in \VV_n$, if $x$ is $n$-successful then
\begin{equation}
 \TT ( i^{-1}(x), \ep_{n,1}) \geq (\log 
\ep_{n,1}
)^2 \Big( \frac{a}{\pi}
 -\frac{2}{
\sqrt{\log n} } \Big).
\label{goalb}
\end{equation}
\end{lemma}

\noindent
{\bf Proof of Lemma~\ref{lem1}:} Recall that if $x$ is
$n$-successful then $\TT ( i^{-1}(x), \ep_{n,1})>
\sum_{j=0}^{n_n} \tau^{(j)}$. Hence, using (\ref{c2.10a}) with
$N=n_n=3 a n^2 \log n$, $\de_n=\pi/(a \sqrt{\log n})$, $R=\ep_{n,n}$,
and $r=
\ep_{n,n-1}$ so that
$\log (R/r)=3\log n$
and $R>r^{0.8}$, we see that for some $C>0$ that is independent of $n$,
\beaa
P_x &:=&\PPP^{x_0}\(\TT ( i^{-1}(x), \ep_{n,1})\leq (
\frac{a}{\pi}-\frac{2}{\sqrt{\log n}})(\log 
\ep_{n,1})^2
\,, \, x \, \mbox{is $n$-successful}\) \nn\\
&\leq &\PPP^{x_0}\(\sum_{j=0}^{N}\tau^{(j)}\leq (
\frac{a}{\pi}-\frac{1}{\sqrt{\log n}})(3n\log n)^2 \)\nn\\
&\leq &\PPP^{x_0}\(\frac{1}{N}
\sum_{j=0}^{N}\tau^{(j)}\leq
(1-\de_n) \frac{\log (R/r)}{\pi} \)
\leq  e^{-Cn^2}.
\eeaa
Consequently, the sum of $P_x$ over all $x\in\VV_n$ and
then over all $n$ is finite, and the
Borel-Cantelli lemma then completes the proof of Lemma~\ref{lem1}.\qed

Taking $\VV_n=\{x_{n,k} : k=1,\ldots, M_n\}$, and the
subsequence
$n(j )= j (\log j)^3$, it follows
from \req{p-new2}, \req{goalb} and the Borel-Cantelli lemma that a.s.
\beq{p-new3}
\CC_{\ep_{n(j ),1}} \geq (\log 
\ep_{n(j),1})^2 \Big(
\frac{a}{\pi}-\frac{2}{\sqrt{\log {n(j )}}} \Big) \,,
\end{equation}
for all $j$ large enough. Since $\ep \mapsto \CC_\ep$ is
monotone non-decreasing,
it follows that for any $\ep_{n(j +1),1}\leq \ep\leq \ep_{n(j ),1}$
\[
\frac{\CC_\ep}{(\log \ep)^2} \geq
\frac{\CC_{\ep_{n(j +1),1}}}{(\log \ep_{n(j ),1})^2}
\]
Observing that $(\log 
\ep_{n(j +1),1})/(\log \ep_{n(j ),1})\to 1$ as
$j \to \infty$,
we thus see that \req{p.1} is an immediate consequence of \req{p-new3}.
\qed

\noindent {\bf Remark:} We note for use in Section \ref{RKproof} that
essentially the same proof
shows that for any $\widehat{a}<2$, almost surely, 
\beq{jpatch1}
\sup_{x\in n(j )^{ -4}S}\TT(x,\ep_{n(j ),1}) \geq (\log \ep_{n(j ),1})^2 \Big(
\frac{\widehat{a}}{\pi}-\frac{2}{\sqrt{\log {n(j )}}} \Big) \,,
\end{equation} 
for all $j$ large enough. To see this we need only prove \req{p-new2}
with the sum now going over $j'$ such  that $x_{ n,j'}\in n^{ -4}S$.
This has the effect of
replacing $M_n$ by $n^{ -4}$ times its previous value. Clearly
(\ref{pm3.7}) still holds,
with perhaps a different $\ga>0$. Also, we now have only
$(\ep_{n,n}/\ep_n)^2 \leq C_2 M_n n^{-42}$, but this is enough to establish
\req{jpatch2}. The rest of the proof follows as before.

\section{Proof of the lattice torus covering time conjecture}\label{sec-etc}

To establish Theorem \ref{theo-2lattice} it suffices to prove that
for any $\de>0$
\begin{equation}
\lim _{n\to \ff}\PPP\({\TT_n\over (n\log n)^2}\geq
\frac{4}{\pi} -\de\)=1\label{kr0}
\end{equation}
since the complementary upper bound on $\TT_n$
is already contained in
\cite[Corollary 25, Chapter 7]{Ald-F}
(see also the references therein).
Our approach is to use Theorem \ref{theo-1p} together with the strong
approximation results of
\cite{Ei}
and
\cite{KMT}.

Fix $\ga>0$ and let 
$\eps_n=2 n^{\ga-1}$.
Then by Theorem \ref{theo-1p} for all $n\geq N_0$ with some
$N_0=N_0(\ga,\de)<\ff$
\begin{equation}
\PPP\(\CC_{\eps_n}> {2(1-\ga-\de)^2\over\pi}(\log n)^2 \)\geq 1-\de
.\label{kr1}
\end{equation}

By Einmahl's~\cite[Theorem 1]{Ei} multidimensional
extension of the Koml\'{o}s-Major-Tusn\'{a}dy~\cite{KMT}
strong approximation theorem,
 we may, for each $n$, construct $\{S_k\}$ and
$\{W_t\}$ on the same probability space so that a.s.
for some $n_0=n_0(\om)<\ff$,
$$
\max_{k\leq 4n^2(\log n )^2}
|W_{k}-{\sqrt{2}} S_k|\leq    n^{\ga}/
6,
\hspace{.2in}\forall n\geq n_0
$$

Hence, dividing by $\sqrt{2}n$ we have
$$ \max_{k\leq 4n^2(\log n)^2 }
|{W_{k}\over \sqrt{2}n}-{S_k\over n}|\leq   \eps_n/2,\hspace{.2in}\forall n\geq
n_0
$$
or, using Brownian scaling, we have
\begin{equation}
 \PPP\(\max_{k\leq 4n^2(\log n)^2 }
|W_{k/2n^2}-{S_k\over n}|\geq  \eps_n/2\)\leq \de
\label{kr2}
\end{equation}
for all $n\geq N'_0$ with some $N'_0=N_0'(\ga,\de)<\ff$.

Now, by (\ref{kr1}) we see that with probability at least $1-\de$ some
disc $\DT (x,\eps_n)\subseteq \Bbb{T}^2$ is completely missed by
\[\lc W_{k/2n^2}\,\bbmod \,\Z^2\,;\,k\leq
{4(1-\ga-\de)^2 \over\pi}n^2 (\log n)^2\rc,
\]
hence by (\ref{kr2}) with probability
at least $1-2\de$
 we have that
\[\lc{S_k\over n}\,\,\bbmod \,\Z^2\,;\,k\leq
{4(1-\ga-\de)^2\over\pi}n^2( \log n)^2\rc\]
avoids some disc of radius
$\eps_n/2=n^{\ga-1}$. Thus, the probability that
\[\lc S_k\bmod
n\Z^2\,;\,k\leq {4(1-\ga-\de)^2\over\pi}n^2 (\log n)^2 \rc\]
avoids some disc of radius
$n^{\ga}$ is at least $1-2\de$, which implies (\ref{kr0}).

\section{Proof of the Kesten-R\'{e}v\'{e}sz conjecture}\label{RKproof}

Let $D_r=D(0,r)\cap \Z^2$ denote the disc of radius $r$ in $\Z^2$ and
define its boundary
\[\partial D_r=\{z\notin D_r\,\Big|\,|z-y|=1\, \mbox{\rm for some } y\in
D_r \}.\]

Let $\phi_n=(\log n)^2/\log\log n$ and let $\NN_n$ denote the number of
excursions
in $\Z^2$ from $\partial D_{2n}$ to
$\partial D_{n (\log n)^3}$
after first hitting $\partial D_{n (\log n)^3}$, that is
needed to
cover $D_{n}$.
By \cite[Theorem 1.1]{Lawler}, it suffices to show that
$$
\limsup_{n\to\infty} \PPP(\log T_n\leq t(\log n)^2) \le e^{-4/t},
$$
and by \cite[Equation (7), page 196]{Lawler}, this is a direct
consequence of the next lemma.
\begin{lemma}\label{lem-RKP}
\begin{equation}
\liminf_{n\rar\ff}{\NN_n\over  \phi_n}
\geq
{2\over
3}\hspace{.1in}\mbox{in probability.}
\label{14.11}
\end{equation}
\end{lemma}

\noindent
{\bf Remark:} Though not needed for our proof of Theorem \ref{theo-4a},
it is not hard to modify  the proof of Lemma \ref{lem-RKP} so as to
show that $\NN_n/\phi_n \to {2\over 3}$ in probability.

Let $K(z,u)$ denote the Poisson kernel for the
annular region $\AA_r :=\{ z : r < |z| < 1/2 \}$,
such  that for any continuous function
$g \geq 0$ on $\partial \AA_r$, we have
$$
\E^z(g(W_\theta)) = \int_{\partial \AA_r} g(u) K(z,u) du \,,
$$
where $\theta:= \inf \{ t \geq 0 : W_t \in \partial \AA_r \}$,
and $W_t$ is a planar Brownian motion, starting at
$W_0=z \in \AA_r$. A preliminary step in
proving Lemma \ref{lem-RKP} is the following estimate about
$K(z,u)$ when $|z| \gg r =|u|$.
\begin{lemma}
\label{lem-harm}
There exists finite $c>2$
such that if $cr \leq |z| <1/(2c)$, then
\beq{har-anul}
\sup_{\{u: |u|=r\}}\, K(z,u) \leq \Big(1 + \frac{40 r \log (2r)}
{|z| \log (2|z|)}
\Big)
\inf_{\{u: |u|=r\}} \,K(z,u)
\end{equation}
\end{lemma}

\noindent
{\bf Proof of Lemma~\ref{lem-harm}:} The series expansion
$$
P_A(x,u)=c_0(x) + \sum_{m=1}^\infty c_m(x) Z_m(x,\frac{u}{|u|})
$$
is provided in
\cite[10.11--10.13, Page 191]{axler} for the Poisson kernel
$P_A(\cdot,\cdot)$ in the region $A=\{ x : r_0 < |x| < 1 \}$,
at its inner boundary $|u|=r_0$, where
$$
c_m(x) = |x|^{-m} \lc \frac{r_0}{|x|} \rc^m \frac{1-|x|^{2m}}{1-(r_0)^{2m}},
\qquad m \geq 1,
$$
and the ``zonal harmonic'' functions
$$
Z_m(x,e^{i\phi})=2 |x|^{m} \cos( m ( \mbox{Arg}(x) - \phi ) )
$$
are given in \cite[5.9 and 5.18]{axler}. Note that for any $x \in A$
\beqn{ofer-1}
|P_A(x,u)-c_0(x)| &\leq& \sum_{m=1}^\infty c_m(x) |Z_m(x,\frac{u}{|u|})| \nn \\
                  &\leq& 2 \sum_{m=1}^\infty \(\frac{r_0}{|x|}\)^m
                  = \frac{2 r_0}{|x|-r_0} \;.
\eeqn
The function $c_0(x)=\log(1/|x|)/\log(1/r_0)$ is the
harmonic function in $A$ corresponding to the
boundary condition $\one_{|x|=r_0}$.
By Brownian scaling $K(z,u)=P_A(2z,2u)$ for $r_0=2r$. Hence, it follows
from \req{ofer-1} and the value of $c_0(\cdot)$,
that for all $2r \leq |z| < 1/2$,
$$
\sup_{\{u: |u|=r\}}\, K(z,u)
\leq \(1 +  \frac{8 f(r)}{f(|z|)-4 f(r)} \)
\inf_{\{u: |u|=r\}}\, K(z,u)\,,
$$
where $f(t):=t \log (1/(2t))$. The proof is complete by
noting that $f(t) \geq 5 f(r)$ for all $cr \leq t \leq 1/(2c)$ provided
$c$ is large enough ($c=10$ suffices).
\qed

With $\Bbb{T}^2=(-1/2,1/2]^2$,
our application of Lemma \ref{lem-harm} is via the following estimate.
\begin{lemma}\label{rn-bd}
Assume $W_0=X_0=\beta$ with $|\beta|=R \in (r,1/2)$, and
let $\tau_r :=\inf\{ t \geq 0: |W_t| = r\}$.
There exists finite $c>2$, such that if $cr \leq R <1/(2c)$, then
the law of $W_{\tau_r}$ is absolutely continuous with respect to
the law of $X_{\tau_r}$,
with  Radon-Nikodym derivative $h_r(\beta,\cdot)$ such that
\beq{ofer-3}
\sup_{|\beta|=R, |\alpha|=r} h_r(\beta,\alpha)
\leq 1 + \frac{40 r \log (2r)}{R \log (2R)} \,.
\end{equation}
\end{lemma}

\noindent
{\bf Proof of Lemma \ref{rn-bd}:}
Recall that the exit time $\theta$ from the annular region
$\AA_r$ is such that $\theta \leq \tau_r$, with equality iff
the path exits $\AA_r$ via its inner boundary $\partial D(0,r)$.
Moreover, with $X_0=W_0=z \in \AA_r$, the path
$\{ X_t : 0 \leq t \leq \theta\}$ is identical in law to
$\{ W_t : 0 \leq t \leq \theta\}$. Let $L$ denote the
number of excursions of $\om_t$ between
$\partial D(0,R)$ and $\partial D(0,1/2)$ completed by time $\tau_r$.
For each $k \geq 0$, let $\mu_k(\beta,\cdot)$
denote the hitting (probability) measure of $\partial D(0,R)$ induced by
$W_t$ upon completing $k$ such excursions, conditional upon $L \geq k$.
Let $\nu_k(\beta,\cdot)$ denote the corresponding hitting measure induced by
the process $X_t$.
Note that $L$ has a Geometric($p$) law, where $p < 1$ is the same
for both processes $X_t$ and $W_t$ and is independent upon the initial
condition $z \in \partial D(0,R)$.  Consequently, for any Borel set
$B \subset \partial D(0,r)$,
\beaa
&& \PPP^\beta(W_{\tau_r} \in B) =
\sum_{k=0}^\infty \PPP^\bb ( W_{\tau_r} \in B,\, L=k)
\\
&=& \sum_{k=0}^\infty p^k \int_{\partial D(0,R)} \mu_k(\bb,dz)
\int_B K(z,u) du \leq \frac{1}{1-p} \int_B [ \sup_{|z|=R} K(z,u) ] du
\,,
\eeaa
where $K(z,u)$ is the Poisson kernel for $W_t$ and the region $\AA_r$.
Similarly,
\beaa
\PPP^\bb(X_{\tau_r} \in B) &=&
\sum_{k=0}^\infty p^k \int_{\partial D(0,R)} \nu_k(\bb,dz)
\int_B K(z,u) du
\\
&\geq& \frac{1}{1-p} \int_B [ \inf_{|z|=R} K(z,u) ] du
\,.
\eeaa
Hence, for any $B \subset \partial D(0,r)$,
$$
\PPP^\bb(W_{\tau_r} \in B) \leq P^\bb(X_{\tau_r} \in B)
     \frac{\sup_{|z|=R, |u|=r} \; K(z,u)}{\inf_{|z|=R, |u|=r} \; K(z,u)}
\,,
$$
implying that $W_{\tau_r}$ is absolutely continuous with respect to
$X_{\tau_r}$, and by \req{har-anul}
 the Radon-Nikodym derivative $h_r(\beta,\cdot)$ clearly
satisfies \req{ofer-3}.
\qed

\noindent
{\bf Proof of Lemma~\ref{lem-RKP}:}
For any $K\subseteq \Bbb{T}^2$ let
\[\CC_\ep(K)=\sup_{x\in K}\TT(x,\ep)\]
 be the $\ep$-covering time of $K$. 
Fix $a>0$ and $b \in (0,1)$. Set $r_\ep=a/|\log \ep|^3$.
Taking the isometry $i:\DT(0,1/2) \mapsto D(0,1/2)$ to be the
identity, omitting $i^{-1}$ throughout the proof,
we can find sequences 
$n(j )\uparrow \ff$ and $\ep_{n(j ),1}\downarrow 0$ with
$(\log \ep_{n(j +1),1})/(\log \ep_{n(j),1})\to 1$ such that 
for any $\widehat{a}<2$, almost surely
$$
\frac{\CC_{\ep_{n(j ),1}}(D(0,b r_{ \ep_{n(j+1 ),1}}))} {\left(\log
\ep_{n(j ),1}\right)^2} \geq  \Big(
\frac{\widehat{a}}{\pi}-\frac{2}{\sqrt{\log {n(j )}}} \Big) \,,
$$
for all $j$ large enough. Indeed, this
follows from \req{jpatch1} after noting that
$n( j)^{ -4}S\subseteq D(0,b r_{ \ep_{n(j+1 ),1}})$.
By first interpolating  for $\ep_{n(j+1 ),1}\leq \ep\leq \ep_{n(j ),1}$
using monotonicity and then letting $\widehat{a}\uparrow 2$ we thus have
that almost surely,
\begin{equation}
\lim_{\ep \rar 0}\frac{\CC_{\ep}(D(0,b r_\ep))}
{\left(\log \ep\right)^2}=\frac{2}{\pi}\,.
\label{14.10}
\end{equation}

Fix $1>\ga>0$. For the remainder of this section only we set
$\ep_n=n^{\ga-1}$ and $r_n=r_{\ep_n}$.
Using the notations of Section \ref{sec-hit},
for $x=0$, $r=r_n$ and any $R \in (0,1/2)$, let
$$
\NN'_n(a,R,b) = \max \{ j : \frak{T}_j \leq \CC_{\ep_n}( D(0,b r_n)) \}
$$
denote the number of excursions of the Brownian motion $X_t$
in the torus $\Bbb{T}^2$ from $\pDT (0,r_n) = \partial D (0, r_n)$ to
$\pDT (0,R) = \partial D(0,R)$
up to time $\CC_{\ep_n}( D(0,b r_n))$.
Fixing $\de>0$, let $N_n=(2/3)(1-\ga)^2 (1-2\de) \phi_n$, noting that
$$
\frac{2}{\pi}(1-\de)(\log \ep_n)^2 \geq (1+\de) \frac{N_n}{\pi} \log(R/r) \,,
$$
for all $n \geq n_0(a,R,\de,\ga)$, implying that,
\beaa
\PPP ( \NN'_n(a,R,b) \leq N_n) &\leq&
\PPP\(
\CC_{\ep_n}(D(0,b r_n)) \leq
\frac{2}{\pi}(1-\de)(\log \ep_n)^2 \) \\
&+&
\PPP \(\sum_{j=0}^{N_n} \tau^{(j)}\geq (1+\de){\log(R/r)\over \pi}N_n \)
\eeaa
Hence, by \req{c2.10b} and (\ref{14.10}) it follows that
for any $R<R_1(\de)$, $a>0$ and $b \in (0,1)$,
\begin{equation}
\lim_{n\rar\ff} \PPP\( \NN'_n(a,R,b) \leq N_n \) = 0 .
\label{14.11a}
\end{equation}
Our next task is to show that \req{14.11a} applies for
the excursion counts $\NN_n(a,R,b)$
that correspond to $\NN'_n(a,R,b)$, when $X_t$ is replaced by
the planar Brownian motion $W_t$.
To this end, consider the random vectors
${\bf W}_k :=(W_{\frak{T}_{j-1}+\sigma^{(j)}}, j=1,\ldots,k)$ and
${\bf X}_k :=(X_{\frak{T}_{j-1}+\sigma^{(j)}}, j=1,\ldots,k)$.
Recall that the $j$-th excursion of $X_t$ from
$\pDT (0,r)$ to $\pDT (0,R)$,
starting at $\alpha_j =X_{\frak{T}_{j-1}+\sigma^{(j)}}$
is precisely the isomorphic image of a planar Brownian motion
started at $\alpha_j$, and run till first hitting $\partial D(0,R)$
(and same applies in case of $\alpha_0=X_0=0$).
Thus, by the strong Markov property of both $X_t$ and $W_t$ at
the stopping times $\frak{T}_{0},
\frak{T}_0+\sigma^{(1)},\frak{T}_1,\frak{T}_1+\sigma^{(2)},\ldots$
we see that for every Borel set $B \subset (\partial D(0,r))^k$
$$
\PPP^0({\bf W}_k \in B) = \E^0 (
\prod_{j=0}^{k-1}
h_r(X_{\frak{T}_j}, X_{\frak{T}_j+\sigma^{(j+1)}} )
;\, {\bf X}_k \in B ) \,.
$$
Recall that $|X_{\frak{T}_j}| = R$ and
$|X_{\frak{T}_j+\sigma^{(j+1)}}| = r$ for all $j \geq 0$.
Consequently,
the law of ${\bf W}_k$ is absolutely continuous with respect to
the law of ${\bf X}_k$,
with Radon-Nikodym derivative $h_{k,r}$ such that
$$
\| h_{k,r} \|_\infty \leq
\( \sup_{|\beta|=R, |\alpha|=r} h_r(\beta,\alpha) \)^k \;.
$$
With $r=r_n \to 0$, we thus have by \req{ofer-3} that for
small enough $R>0$ and all $n$ large enough,
\beq{ofer-2}
\| h_{N_n,r_n} \|_\infty \leq \( 1 + \frac{40 r_n \log (2r_n)}{R \log (2R)}
\)^{N_n} \;.
\end{equation}
Since $N_n r_n |\log (2 r_n)| \to 0$, we see that
$ \| h_{N_n,r_n} \|_\infty \to 1$ as $n \to \infty$.
Since $b<1$, and with the $j$-th excursion of $X_t$ from
$\pDT (0,r)$ to $\pDT (0,R)$,
starting at some $\alpha_j =X_{\frak{T}_{j-1}+\sigma^{(j)}}$
being the isomorphic image of a planar Brownian motion
started at $\alpha_j$, and run till first hitting $\partial D(0,R)$,
we get by the strong Markov property of both $X_t$ and $W_t$ that
for any $k$,
$$
\E \( \one_{\NN_n(a,R,b) \leq k} \,|\, \sigma({\bf W}_k) \) =
\E \( \one_{\NN'_n(a,R,b) \leq k}\,|\, \sigma({\bf X}_k) \) \,,
$$
implying that
\beq{ofer-4}
\PPP \( \NN_n(a,R,b)\leq k \) = \E \( h_{k,r_n} ({\bf X}_k) \,,\,
 \NN'_n(a,R,b)\leq k \)
\end{equation}
It thus follows from \req{14.11a}, \req{ofer-2} and \req{ofer-4} that
\beqn{14.11b}
\PPP\( \NN_n(a,R,b)\leq N_n \) &=&
\E\( h_{N_n,r_n}({\bf X}_{N_n} )
 \,,\, \NN'_n(a,R,b)\leq N_n \) \nn\\
&\leq&
 \| h_{N_n,r_n} \|_\infty \PPP\(  \NN'_n(a,R,b) \leq N_n \) \to 0
\eeqn
Setting $R<R_0(\de)$ small enough for \req{14.11b} to apply,
with $a:=2 R(1-\ga)^5$ and $b:=1/(2(1-\ga))$, we next
use strong approximation, as in Section \ref{sec-etc}, to show how
(\ref{14.11}) follows from this. Indeed, with $t_n := \exp((\log n)^3)$,
we may and shall, for each $n$, construct $\{S_k\}$ and
$\{W_t\}$ on the same probability space so that for some $n_0=n_0(\om)<\ff$
$$ \max_{k\leq t_n}
|W_{k}-{\sqrt{2}} S_k|\leq    n^{\ga/2},\hspace{.2in}\forall n\geq
n_0\hspace{.2in} a.s.$$
Hence, multiplying by $\rho_n:= b r_n/(\sqrt{2} n)$ we have
$$ \max_{k\leq t_n}
|\rho_n W_{k}- \rho_n \sqrt{2} S_k|\leq
\eps_n/3,\hspace{.2in}\forall n\geq n_0\hspace{.2in} a.s.$$ or, using
Brownian scaling, we have
\begin{equation}
 \PPP\(
\max_{k\leq t_n}
|W_{k \rho_n^2}- \rho_n \sqrt{2} S_k|\leq  \eps_n/3\)\geq
1-\de\label{14.kr2}
\end{equation}
for all $n\geq N'_0$ with some $N'_0=N_0'(\ga,\de)<\ff$.

Recall that $\PPP(T_n > t_n) \to 0$, see
\cite[Theorem 1.1]{Lawler}, hence
by \req{14.11b}, we see that for all $n$ sufficiently large,
\begin{equation}
\PPP\(\NN_n(a,R,b)>N_n,\, T_n \leq t_n\)\geq 1-\de \,.
\label{14.20}
\end{equation}
Now, by (\ref{14.20}) we have that with
probability at least $1-\de$ some disc
$D(x,\eps_n)\subseteq D(0,b r_n)$ is completely missed by
$\{ W_{k \rho_n^2} \}$ during the first $N_n$ excursions from
$\partial D(0,r_n)$ to $\partial D(0,R)$.
Moreover,  by \req{14.20}, also
$\{ \sqrt{2} \rho_n S_k : k \leq t_n\}$ covers
$\sqrt{2} \rho_n D_n$, hence with probability at least $1-2\de$,
we also have by \req{14.kr2}, that the sequence
$\{ W_{k \rho_n^2} : k \leq t_n\}$  provides
a $(2 \ep_n/3)$-cover of the set $D(0,\sqrt{2} \rho_n n)$.
Our choice of $\rho_n$ guarantees that the latter set is exactly
$D(0,b r_n)$. Consequently, in this case we know that
the $N_n$ excursions mentioned above
are completed by time $\rho_n^{-2} t_n$.
Observe that $b>1/2$ and $b r_n (\log n)^3 = R (1-\ga)$, hence
$(r_n+\ep_n/3) <
\sqrt{2} \rho_n (2n)$ and
$(R-\ep_n/3) >
\sqrt{2} \rho_n n (\log n)^3$, for all $n$ large.
Appealing again to (\ref{14.kr2}) we thus further have that
$\{ \sqrt{2} \rho_n S_k \}$ avoids some disc of radius
$\eps_n/3={1\over 3}n^{\ga-1}$ in $D(0,\sqrt{2} \rho_n n)$
during its first $N_n$ excursions from $\sqrt{2} \rho_n \partial D_{2n}$
to $\sqrt{2} \rho_n \partial D_{n (\log n)^3}$. Thus, the probability
that $\{ S_k \}$ avoids some lattice point
in $D_n$  during its first $N_n=\frac{2}{3}(1-\ga)^2 (1-2\de) \phi_n$
excursions from
$\partial D_{2n}$ to $\partial D_{n (\log n)^3}$ is at least
$1-2\de$. Considering $\de \to 0$, followed by $\ga \to 0$,
we get \req{14.11}.
\qed

\section{First moment estimates}\label{1momest}

We start with analyzing the birth-death Markov chain $\{ Y_l \}$ on
the state space $\{-n,-(n-1),\ldots,-1\}$, starting at $Y_0=-n$,
having both $-n$ and $-1$ as reflecting boundaries (so that
$\PPP(Y_l=-(n-1)|Y_{l-1}=-n)=1$, $\PPP(Y_l=-2|Y_{l-1}=-1)=1$)
and the transition probabilities
\bea
\label{pkdef}
\oo{p}_{k} := \PPP(Y_l=-(k-1) | Y_{l-1}=-k)
&=& 1-\PPP(Y_l=-(k+1) | Y_l=-k) \nonumber \\
&=&
{\log (k+1) \over \log k+\log (k+1)}
\;.
\eea
for $k=2,\ldots, n-1$. Let $\ze=3a>0$ and
\[
\sss
:=\inf \{m\,:\,\sum_{j=1}^m \one_{\{-n\}} (Y_{j}) =\ze n^2\log n\},
\]
denote the number of steps it takes this birth-death Markov chain
to complete $\ze n^2\log n$ excursions from $-(n-1)$ to $-n$.
For each $-n \leq k \leq -2$,
$$
\oo{L}_k = \sum_{l=1}^{\sss} \one_{\{Y_{l-1}=k,Y_l=k+1\}} \;,
$$
denote the number of transitions of $\{ Y_l \}$ from state $k$ to state
$k+1$ up to time $\sss$. (Thus, $ \oo{L}_{-n}=\ze n^2\log n$).
As we show below, fixing $x \in S$,
the law of $\{ N^x_{n,k}\}_{k=2}^{n}$ relevant for the $n$-successful
property, is exactly that of $\{ \oo{L}_{-k} \}_{k =2}^{n}$.
To get a hold on the latter, note that
conditional on $\oo{L}_{-(k+1)}=\ell_{k+1} \geq 0$ we
have the representation
\begin{equation}
\oo{L}_{-k}=\sum_{i=1}^{\ell_{k+1}} Z_i \,,
\label{m1.4}
\end{equation}
where the $Z_i$ are independent identically distributed (geometric) random
variables with
\begin{equation}
\PPP(Z_i=j)= (1-\oo{p}_{k}) \oo{p}_{k}^j \;,
\hspace{.3in}j=0,1,2,\ldots\label{m1.5}
\end{equation}
Consequently, $\{ \oo{L}_k \}_{k=-n}^{-2}$ is a Markov chain on $\Z_+$
with initial condition $\oo{L}_{-n}=\ze n^2\log n$,
and transition probabilities
$\PPP(\oo{L}_{-k}=0|\oo{L}_{-(k+1)}=0)=1$,
\begin{equation}
\PPP\(
\oo{L}_{-k}=\ell\,\big|\,
\oo{L}_{-(k+1)}=\wt{m} \)={\wt{m}-1+\ell \choose \wt{m}-1}
\oo{p}_{k}^{\ell}(1-\oo{p}_{k})^{\wt{m}},
\label{m1.6}
\end{equation}
for $\wt{m} \geq 1$, $\ell \geq 0$ and $k=n-1,\ldots,2$.

Let $n_k=\ze k^2\log k$ for $k=3,\ldots,n-1$ and define
for $2 \leq i < j \leq n$,
\beq{def-h}
h_{i,j}(\ell_j) := \sum_{\stackrel{\ell_{i},\ldots,\ell_{j-1}}{
|\ell_k-n_k| \leq k}}
\prod_{k=i}^{j-1}\PPP \( \oo{L}_{-k}=\ell_{k} \,
\big|\,\oo{L}_{-(k+1)}=\ell_{k+1}\) ,
\end{equation}
where $\ell_n=\ze n^2\log n$ and $\ell_2=0$.
The next lemma is key to estimating the
growth of $h_{i,n}(\ell_n)$ in $n$.
\begin{lemma}\label{recursion}
For some $C=C(\ze) <\infty$ and all $3\leq k\leq n-1$,
$|\ell-n_k| \leq k$, $|\wt{m}-n_{k+1}| \leq  k+1$, $\wt{m} \geq 1$,
\begin{equation}
C^{-1} {k^{-(\ze+1)} \over \sqrt{\log k}} \leq
\PPP\(\oo{L}_{-k}=\ell\,\big|\,\oo{L}_{-(k+1)}=\wt{m} \) \leq C
{k^{-(\ze+1)}\over \sqrt{\log k}}\,.
\label{m1.2}
\end{equation}
\end{lemma}

\noindent
{\bf Proof of Lemma \ref{recursion}}: With
$p_k=1-\oo{p}_k$ and $m=\wt{m}-1 \geq 0$, we see that
\beq{amir-p7.6}
\frac{1-p_k}{p_k} \PPP\(\oo{L}_{-k}=\ell\,\big|\,\oo{L}_{-(k+1)}=\wt{m}\)
= {m+\ell \choose m } p_{k}^{m} (1-p_{k})^{\ell+1} \;.
\end{equation}
The right hand side of \req{amir-p7.6} is merely \cite[(7.6)]{DPRZ4}
for which the bounds of \req{m1.2} are derived in \cite[Lemma 7.2]{DPRZ4}.
To complete the proof, note that
$p_k=1-\oo{p}_k$ is bounded away from $0$ and $1$
(see \req{pkdef}).
\qed

Note that
\[
\inf_{\wt{m} \leq n_3+3}
\PPP\(\oo{L}_{-2}=0\,\big|\,\oo{L}_{-3}=\wt{m}
\) \geq (1-\oo{p}_2)^{n_3+3}>0.
\]
Hence, setting $h_{n,n}(\ell_n)=1$,
it follows from \req{def-h} and \req{m1.2} that
for some $C_1<\infty$,
\beq{rel-bd}
C_1^{-1} \frac{k^{-\ze}}{\sqrt{\log k}}
 \leq \frac{h_{k,n}(\ell_n)}{h_{k+1,n}(\ell_n)}
\leq C_1 \frac{k^{-\ze}}{\sqrt{\log k}}
 \qquad \forall \;\; 2 \leq k \leq n-1\,.
\end{equation}
Applying \req{rel-bd} we conclude also that
for any $\ga>0$ there exists $C_2=C_2(\ga)>0$ such that
for all  $2 \leq l \leq n-1$.
\beq{h-bd}
h_{l,n}(\ell_n) \geq \prod_{k=l}^{n-1} C_1^{-1}
\frac{k^{-\ze}}{\sqrt{\log k}}
\geq C_2^{n-l} \lc \frac{n!}{l!} \rc^{-\ze-\ga}
\end{equation}

Recall that $\ep_k=\ep_1 (k!)^{-3}$ and
$\ep_{n,k} = \rho_n \ep_n (k!)^3$
for $\rho_n=n^{-21}$ and $k=1,\ldots,n$.
Per $n \geq 3$ and $x \in S =[\ep_1,2\ep_1]^2$,
$\RR^x_n$ denotes the time until $X_t$ completes
$\ze n^2\log n$ excursions from
$i^{-1}(\partial D(x,\ep_{n,n-1}))$ to
$i^{-1}(\partial D(x,\ep_{n,n}))$
and $N_{n,k}^x$, $k=2,\ldots,n$,
denote the number of excursions from
$i^{-1}(\partial D(x,\ep_{n,k-1}))$ to
$i^{-1}(\partial D(x,\ep_{n,k}))$
until $\RR^x_n$. A point $x\in S$ is
{\bf $n$-successful} if
$$
N^x_{n,2}=0, \qquad n_k-k \leq N^x_{n,k} \leq n_k+k  \quad
\forall k=3,\ldots,n-1\,.
$$

The next lemma applies \req{rel-bd} to estimate the first moment of
the $n$-successful property.
\begin{lemma}\label{moment1}
For all $n\geq 3, x\in S$ and some $\de_n \to 0$,
independent of $\rho_n$,
\begin{equation}
\bar{q}_n :=\PPP(x\, \mbox{ is $n$-successful}) = (n!)^{-\ze-\de_n} \,.
\label{m1.1}
\end{equation}
\end{lemma}

\noindent
{\bf Proof of Lemma \ref{moment1}}: Observe that
$$
\oo{p}_{k}=
{\log (\ep_{n,k+1}/\ep_{n,k})\over \log (\ep_{n,k+1}/\ep_{n,k-1})}
\,,
$$
is exactly the probability that
the planar Brownian motion $B_t$ starting at
any $z \in \partial D(x,\ep_{n,k})$
will hit $\partial D(x,\ep_{n,k-1})$ prior to
hitting $\partial D(x,\ep_{n,k+1})$, with $(Y_{l-1},Y_l)$ recording
the order of excursions the
Brownian
path makes between the sets
$\{ \partial D(x,\ep_{n,k}), n \geq k \geq 1 \}$.
Note that $0 \notin D(x,\ep_1)$ for $x \in S$,
the above mentioned probabilities are independent of
the starting points of the excursions, and
$\partial D(x,\ep_{n,k}) \subset D(x,\ep_1) \subset
D(0,1/2)$,
for all $k=1,\ldots,n$. Hence, by the strong Markov property
of the Brownian motion $X_t$ on $\Bbb{T}^2$
with respect to the starting times of its first $n_n$
excursions from
$i^{-1}(\partial D(x,\ep_{n,n-1}))$ to
$i^{-1}(\partial D(x,\ep_{n,n}))$,
it follows
that in computing $\bar{q}_n$ of (\ref{m1.1}) we may and shall replace
$X_t$ by the planar Brownian motion
$B_t =i(X_t) $.
It follows from radial symmetry and the
strong Markov property of Brownian motion that
$\bar{q}_n$ is independent of $x\in S$.
By Brownian scaling, $\bar{q}_n$ is also independent of the value of
$\rho_n \leq 1$. Moreover, as already mentioned, fixing $x \in S$,
the law of $\{ N^x_{n,k}\}_{k=2}^{n}$ is exactly that of
$\{ \oo{L}_{-k} \}_{k =2}^{n}$. We thus deduce that
\beq{m1.1a}
\bar{q}_n = \PPP\(  |\oo{L}_{-k}-n_{k}| \leq k\,;\,3 \leq k\leq
n-1\,;\,\oo{L}_{-2}=0\) = h_{2,n}(\ell_n)
\end{equation}
Since $n^{-1} \log n! \to \infty$ and for some $\eta_n \to 0$
$$
\prod_{k=2}^n \log (k) = (n!)^{\eta_n} \;,
$$
we see that the estimate (\ref{m1.1}) on $\bar{q}_n$ is a
direct consequence of the bound (\ref{rel-bd}).
\qed

In Section \ref{2momest} we control the second moment of the
$n$-successful property. To do this, we need to
consider excursions between disks centered at $x \in S$ as well as
those between disks centered at $y \in S$, $y \neq x$.  The radial
symmetry we used in proving Lemma \ref{moment1} is hence lost. The next
Lemma shows that, in terms of the number of excursions, not much is lost
when we
condition on a certain $\si$-algebra ${\mathcal G}_l^y$ which
contains more information
than just the number of excursions in the previous level.
To define $\GG_l^y$, let
$\tau_0=0$ and for $i=1,2,\ldots$ let
\begin{eqnarray*}
\tau_{2i-1} & = & \inf \{ t \geq \tau_{2i-2}  :\; X_t \in
i^{-1}(\partial D(y,\ep_{n,l-1}))
\} \\
\tau_{2i} & = & \inf \{ t \ge \tau_{2i-1} :\; X_t \in
i^{-1}(\partial D(y,\ep_{n,l}))
\}.
\end{eqnarray*}
Thus, $N_{n,l}^y=\max \{i\,:\,\tau_{2i}\leq \RR_n^y\}$.  For each
$j=1,2,\ldots,N_{n,l}^y$
let
\[e^{(j)}=\{ X_{\tau_{2j-2}+t} : 0\leq t \leq \tau_{2j-1}-\tau_{2j-2} \}\]
be the $j$-th excursion from
$i^{-1}(\partial D(y,\ep_{n,l}))$ to
$i^{-1}(\partial D(y,\ep_{n,l-1}))$
(but note that for $j=1$ we do begin at $t=0$). Finally, let
\[
e^{(N_{n,l}^y+1)}= \{ X_{\tau_{2N_{n,l}^y}+t} : t\geq 0 \}.
\]
We let $J_l:=\{ l-1,\ldots,2 \}$ and take
$\GG_l^y$ to be the $\sigma$-algebra generated by the excursions
$e^{(1)}, \ldots,e^{(N_{n,l}^y)}$, $e^{(N_{n,l}^y+1)}$.


\begin{lemma}\label{recursionends}
For some $C_0<\infty$, any $3 \leq l \leq n$,
$|m_l-n_l| \leq l$ and all $y \in S$,
\bea
&& \PPP(N^y_{n,k} = m_k ; k \in J_l \,
|\,N_{n,l}^y=m_l, {\mathcal G}_l^y) \nonumber \\
&&\leq
(1+C_0 l^{-1} \log l)  \prod_{k=2}^{l-1} \PPP\(\oo{L}_{-k}=m_{k} \,|\,
\oo{L}_{-(k+1)}=m_{k+1} \)
\label{m1.2e}
\eea
\end{lemma}
The key to the proof of Lemma \ref{recursionends} is to
demonstrate that the number of Brownian excursions involving
concentric disks of radii $\ep_{n,k}$, $k \in J_l$ prior to
first exiting the disk of radius $\ep_{n,l}$ is almost
independent of the initial and final points of the
overall excursion between the $\ep_{n,l-1}$ and $\ep_{n,l}$ disks.
The next lemma provides uniform estimates sufficient for this task.

\begin{lemma}\label{lemprobends}
Consider a Brownian path $B_\cdot$
starting at $z\in \partial D(y,\ep_{n,l-1})$,
for some $3 \leq l \leq n$. Let
$\bar{\tau}=\inf \{t>0\, :\,B_t\notin D(y,\ep_{n,l})\}$ and
$Z_k$, $k \in J_l$,
denote the number of excursions of the path  from
$\partial D(y,\ep_{n,k-1})$ to  $\partial D(y,\ep_{n,k})$, prior to
$\bar{\tau}$.
Then, there exists a universal constant $c<\infty$, such that for
all $\{ m_k : k \in J_l \}$,
uniformly in $v\in \partial D(y,\ep_{n,l})$ and $y$,
\begin{equation}\qquad
\PPP^z (Z_k = m_k , k \in J_l \,\big |\,B_{\bar{\tau}}=v)
\leq (1 + c l^{-3}) \PPP^z (Z_k = m_k ,\, k \in J_l) \,.
\label{m1.5eu}
\end{equation}
\end{lemma}

\noindent
{\bf Proof of Lemma \ref{lemprobends}:} This is essentially
\cite[Lemma 7.4]{DPRZ4}. The only difference is that here
we use the sequence of radii $\ep_{n,k}$, for $k=l,l-1,l-2,\ldots,2$,
whereas \cite{DPRZ4} uses the radii $\ep_k$, for $k=l-1,l,l+1,\ldots,n$.
The proof of \cite[Lemma 7.4]{DPRZ4}
involves only the ratio $\ep_l/\ep_{l-1}=l^{-3}$ between the two
exterior disks and the fact that the probability
$p_l$ of reaching the next disk (of radius $\ep_{l+1}$ there),
is uniformly bounded away from $1$. The ratio of the
two exterior disks here is $\ep_{n,l-1}/\ep_{n,l}=l^{-3}$ which is the
same as in \cite{DPRZ4}, whereas $p_l$ is replaced here
by $\oo{p}_{l-1}$, which is also uniformly bounded away from $1$.
\qed

\noindent
{\bf Proof of Lemma \ref{recursionends}}:
Fixing $3 \leq l \leq n$ and
$y \in S$, let $Z^{(j)}_{k}$, $k \in J_l$ denote the
number of excursions from
$i^{-1}(\partial D(y,\ep_{n,k-1}))$ to
$i^{-1}(\partial D(y,\ep_{n,k}))$
during the $j$-th excursion of the path $X_t$ from
$i^{-1}(\partial D(y,\ep_{n,l-1}))$ to
$i^{-1}(\partial D(y,\ep_{n,l}))$.
If $m_l=0$ then
the probabilities in
both sides of \req{m1.2e} are zero unless
$m_k=0$ for all $k \in J_l$, in which case they are both one,
so the lemma trivially applies when $m_l=0$.
Considering hereafter $m_l>0$, since
$0 \notin
i^{-1} (D(y,\ep_1))$ we have that conditioned upon
$\{N^y_{n,l}=m_l\}$,
\begin{equation}
N^y_{n,k}=\sum_{j=1}^{m_l} Z^{(j)}_{k}\;\;\;\;\;\; k \in J_l\,.
\label{m1.4e}
\end{equation}
Conditioned upon $\GG^y_l$, the random vectors
$\{ Z^{(j)}_{k}, k \in J_l\}$ are independent for $j=1,2,\ldots,m_l$.
Moreover, with $X_t$
being the isomorphic image of a planar Brownian motion $B_t$
within $D(y,\ep_{n,l})$, we see that
$\{ Z^{(j)}_{k}, k \in J_l\}$ then has the conditional law of
$\{ Z_{k}, k \in J_l\}$ of Lemma \ref{lemprobends}
for some random $z_j \in \partial D(y,\ep_{n,l-1})$ and
$v_j \in \partial D(y,\ep_{n,l})$, both measurable on ${\mathcal
G}^y_l$
(as $z_j$
corresponds to
the final point of $e^{(j)}$, the $j$-th excursion from
$i^{-1}(\partial D(y,\ep_{n,l}))$ to
$i^{-1}(\partial D(y,\ep_{n,l-1}))$
and $v_j$
corresponds to
the initial point of the $(j+1)$-st such excursion $e^{(j+1)}$).
Let $\PP_l$ denote the finite set of
all partitions
$\{ m_k^{(j)}, k \in J_l, j=1,\ldots,m_l :
m_k=\sum^{m_l}_{j=1} m^{(j)}_k, k \in J_l \}$.
Then, by the uniform upper bound of (\ref{m1.5eu})
and radial symmetry,
\beaa
&&
\PPP(N^y_{n,k}  =
m_k,\,k \in J_l\big|\,N^y_{n,l}=m_l, {\mathcal G}_l^y)
\\
&=&
\sum_{\PP_l}\prod_{j=1}^{m_l}\PPP^{z_j} (Z_k = m^{(j)}_k ,\, k \in
J_l \,\big
|\,B_{\bar{\tau}}=v_j)
\\
&\leq & \sum_{\PP_l}\prod_{j=1}^{m_l} \,(1 + c l^{-3})\PPP^{z_j}(Z_k =
m^{(j)}_k ,\, k \in
J_l )\\
& =&(1 + cl^{-3} )^{m_l}
\PPP\(N^y_{n,k}=m_k,\,k \in J_l\big|\,N^y_{n,l}=m_l\) \,.
\eeaa
Since $m_{l} \leq c_1 l^2 \log l$ we thus get the bound
(\ref{m1.2e}) by the representation used in the proof of Lemma
\ref{moment1}. \qed

\section{Second moment estimates}\label{2momest}

Recall that $N_{n,k}^x$ for $x\in S$, $2 \leq k \leq n$, denotes
the number of excursions from
$i^{-1}(\partial D(x,\ep_{n,k-1}))$ to
$i^{-1}(\partial D(x,\ep_{n,k}))$
prior to $\RR_n^x$. With $n_k=\ze
k^2 \log k$ we shall write $N\stackrel{k}{\sim} n_{k}$ if
$|N-n_k| \leq k$ for $3\leq k\leq n-1$ and $N=0$ when $k=2$.
Relying upon the first moment
estimates of Lemmas \ref{moment1} and \ref{recursionends}, we next
bound the second moment of the $n$-successful property.

\begin{lemma}\label{moment2}
For any $\gamma> 0$ we can find $C=C(\gamma)<\ff$ such that
for all $x,y\in S$,
\begin{equation}
\PPP\(x\, \mbox{ and } \, y \, \mbox { are $n$-successful } \)
\leq \bar{q}_n^2 \;n^{3+5\ze} C^{n-l} \ffrac{n!}{l!}^{\ze+\ga} \;,
\label{m2.2}
\end{equation}
where $l=\max\{k \leq n: |x-y| \geq 2 \ep_{n,k} \} \vee 1$
and $\bar{q}_n:=\PPP(x\, \mbox{ is $n$-successful } )$.
Furthermore, if $|x-y| \geq 2 \ep_{n,n}$ then for some $C_0<\infty$,
\begin{equation}
\PPP\(x\, \mbox{ and } \, y \, \mbox { are $n$-successful } \)
\leq (1+C_0 n^{-1} \log n)\bar{q}_n^2  \;.
\label{m2.200}
\end{equation}
\end{lemma}

\noindent
{\bf Proof of Lemma \ref{moment2}}:
Fixing $x,y \in S$, suppose
$2 \ep_{n,l+1} > |x-y| \geq 2 \ep_{n,l}$ for some $n-1 \geq l \geq 3$.
Since $\ep_{n,l+2} - \ep_{n,l} \geq 2 \ep_{n,l+1}$, it is easy to see that
$i^{-1}(D(y,\ep_{n,l})) \cap i^{-1}(\partial D(x,\ep_{n,k})) =\emptyset$
for all $k \neq l+1$.
Replacing hereafter $l$ by $l \wedge (n-3)$, it is easy
to see that for $k \neq l+1$, $k \neq l+2$,
the events $\{ N_{n,k}^x \stackrel{k}{\sim} n_{k} \}$ are measurable
on the $\sigma$-algebra $\GG^y_{l}$
defined above Lemma \ref{recursionends}.
With $J_l:=\{ l-1,\ldots,2 \}$ and $I_l:=\{ 2,\ldots,l,l+3,\ldots,n-1 \}$,
we note that
\begin{equation}
\{ x, y \, \mbox{ are $n$-successful} \}\subset
\{ N_{n,k}^x \stackrel{k}{\sim} n_{k},\; k \in I_l \}\bigcap \{
N_{n,k}^y
\stackrel{k}{\sim} n_{k},\; k \in J_{l+1}\}\,.\nn
\end{equation}
Let $\MM(I_l):= \{m_2,\ldots,m_{n-1} : m_k\stackrel{k}{\sim} n_{k}, k
\in I_l \}$
(note that the range of $m_{l+1}\,,\,m_{l+2}$ is unrestricted), and
$\MM(J_l):=
\{m_{2},\ldots,m_{l-1} : m_k\stackrel{k}{\sim} n_{k}, k\in J_l \}$.
Applying (\ref{m1.2e}),
we have that for some universal constant $C_3<\infty$,
\beqn{2mom}
&& \PPP\(x\, \mbox{ and } \, y \, \mbox { are $n$-successful} \) \nn \\
&\leq\!\!\!&
\sum_{\MM(J_{l+1})}
\E \left[ \PPP ( N_{n,k}^y=m_k,\, k \in J_{l}\, \big| \,
N^y_{n,l} = m_{l} , {\mathcal G}^y_{l} ) \,;
N_{n,k}^x \stackrel{k}{\sim} n_{k},\, k \in I_l \right] \nn \\
&\leq\!\!\!& C_3
\PPP \( N_{n,k}^x \stackrel{k}{\sim} n_{k},\; k \in I_l \)
\sum_{|m_l -n_l| \leq l} h_{2,l}(m_l)
\nn \\
\eeqn
Since,
\beaa
&&
 \sum_{m_{l+1},\,m_{l+2}}
 \prod_{k=l}^{l+2}
\PPP \(\oo{L}_{-k}=m_{k} \, \big|\,
\oo{L}_{-(k+1)}=m_{k+1}\) \\
&&
\hspace*{3cm}
=
 \PPP \( \oo{L}_{-l}=m_{l} \, \big|\,\oo{L}_{-(l+3)}=m_{l+3}\)\leq 1,
\eeaa
taking $m_n=\ze n^2 \log n$, we have by the representation (\ref{m1.1a})
of Lemma \ref{moment1}, that
\begin{eqnarray}
\PPP \( N_{n,k}^x \stackrel{k}{\sim} n_{k},\; k \in I_l \)
&=&
 \sum_{\MM(I_l)}
\prod_{k=2}^{n-1} \PPP \( \oo{L}_{-k}=m_{k} \, \big|\,
\oo{L}_{-(k+1)} =m_{k+1}\)\nn \\
&\leq& h_{l+3,n}(m_n) \sum_{|m_l -n_l| \leq l} h_{2,l}(m_l)
\label{expand}
\end{eqnarray}
(as mentioned, the sum over ${\MM(I_l)}$ involves the unrestricted
$m_{l+1}$ and $m_{l+2}$).
Combining (\ref{2mom}) and (\ref{expand}), we have
\beq{expandco}
\PPP\(x\, \mbox{ and } \, y \, \mbox { are $n$-successful}
\)
\leq C_3 h_{l+3,n} (m_n)
\Big[\sum_{|m_l -n_l| \leq l} h_{2,l}(m_l) \Big]^2
\end{equation}
By (\ref{m1.1a})
and the bounds of Lemma \ref{recursion} we have the inequalities,
\begin{eqnarray}
\bar{q}_n = h_{2,n}(m_n)
&\geq&
h_{l,n}(m_n) \, \inf_{|m_{l}-n_{l}| \leq l} \, h_{2,l}(m_l) \nn\\
&\geq& h_{l,n}(m_n) C^{-2} \, \sup_{|m_{l}-n_{l}| \leq l} h_{2,l}(m_l)\, \nn\\
&\geq&  h_{l,n}(m_n) C^{-2} (2l+1)^{-1} \sum_{|m_l -n_l| \leq l} h_{2,l}(m_l)
\label{1condmom}
\end{eqnarray}
Combining (\ref{expandco}) and (\ref{1condmom}), we see that
for some universal constant $C_4<\infty$,
$$
\PPP\(x\, \mbox{ and } \, y \, \mbox { are $n$-successful} \)
\leq
C_4 n^{2} \bar{q}_n^2 \frac{h_{l+3,n}(m_n)}{h_{l,n}(m_n)^2} \;.
$$
By \req{rel-bd},
$h_{l+3,n}(m_n)/h_{l,n}(m_n) \leq C_5 n^{3\ze+1}$
for some $C_5 <\infty$ and all $l \leq n-3$.
Thus, we get (\ref{m2.2}) via the bound \req{h-bd} on $h_{l,n}(m_n)$,
with the extra $n^{2 \ze}$ factor coming from the use of
$l \wedge (n-3)$ throughout the above proof.
It also follows from \req{h-bd} and \req{m1.1a} that
when $2 \ep_{n,3} > |x-y|$, the trivial bound
$\PPP\(x\, \mbox{ and } \, y \, \mbox { are $n$-successful} \) \leq
\bar{q}_n$ already implies \req{m2.2}.

Suppose next that $|x-y| \geq 2 \ep_{n,n}$, in which case \req{m2.2}
is contained in the sharper bound (\ref{m2.200}). To prove the latter,
note that if $|x-y| \geq 2 \ep_{n,n}$, then the event
$\{x\,\mbox { is $n$-successful }\}$ is ${\mathcal G}^y_{n}$
measurable, hence
\beaa
&&
\PPP\(x\, \mbox{ and } \, y \, \mbox { are $n$-successful } \)
\label{m2.200a}\\
&& =\E\( \lc \PPP (y\,\mbox { is $n$-successful }\, \big| \,
 {\mathcal G}^y_{n} )\rc\,,\,x\,\mbox { is $n$-successful }\)  \nn
\\
&& =\E\( \lc \PPP (\,
N_{n,k}^y \stackrel{k}{\sim} n_{k},\, k \in J_n \,
\big| \, N^y_{n,n}=m_n,\,
 {\mathcal G}^y_{n} )\rc\,,\,x\,\mbox { is $n$-successful }\)  \nn
\,,
\eeaa
and (\ref{m2.200}) follows from Lemma \ref{recursionends}.
\qed

\section{The $\ep$-covering time of a compact Riemannian manifold}
\label{sec-manifold}

Let $M$ be a
smooth,
compact, connected two-dimensional,
Riemannian manifold without boundary.
Let $\{X_t\}_{t \ge 0}$ denote Brownian motion on $M$ starting at
some non-random $x_0 \in M$.
The process $\{X_t\}_{t \ge 0}$ is a
symmetric, strong Markov
process with reference measure given by the Riemannian measure $dA$ and
infinitesimal generator $1/2$ the Laplace-Beltrami operator
$\Delta_M$. We use
$d(x,y)$ to denote the Riemannian distance between $x,y\in M$.
With this notion of distance we can take over the definitions used for the
plane and the flat torus: $D_M(x,r)$ denotes the
open disc in $M$ of radius $r$ centered at $x$. For $x$ in
$M$  we have the $\ep$-hitting time
\[\TT(x,\eps)=\inf \{t>0\,|\,X_t\in D_M(x,\ep)\}.\]
Then
\[\CC_\ep=\sup_{x \in M}\TT(x,\eps) \]
is the $\ep$-covering time of $M$.

\noindent
{\bf Proof of Theorem \ref{theo-m1p}:}
If $g$ denotes the Riemannian metric for $M$,
let $M'$ denote the Riemannian
manifold obtained by changing the Riemannian metric for $M$ to $g'=g/A$,  so
that the area of $M'$ is $1$. Since $\Delta_{M'} = \frac{1}{A} \Delta_M$,
it follows that $X'_t=X_{t/A}$ is the Brownian motion
on $M'$. With $\CC'_{\ep'}$ denoting the $\ep'$-covering time of $M'$, we
see that $\CC_\ep$ has the same law as $A \CC'_{\ep/\sqrt{A}}$.
Consequently, it suffices to prove the theorem only for manifolds
of area $A=1$, which we assume hereafter.
Then, the statement and
proof of Lemma \ref{lem-hit} applies for any fixed $x \in M$,
upon replacing
$\DT (x,\cdot)$ by $D_M(x,\cdot)$.

Our assumptions about $M$ imply the existence for some $\xi>0$
of a smooth isothermal coordinate system in each disc $D_M(u,\xi)$,
$u \in M$ (c.f. for example \cite[Page 386 and Addendum 1]{spivak}).
This implies that with respect to such coordinates, the
Laplace-Beltrami operator $\Delta_M$ is given on $D_M(u,\xi)$ by
$a(z) (\partial_1^2+\partial_2^2)$ for some smooth, scalar
function $a:M \to (0,\infty)$, with $a(z)=a_u(z)$ possibly
depending
on $u$.
Moreover, for each $u \in M$ and $\de>0$, upon choosing
$\xi=\xi(u,\de)>0$ small enough, we may
after translation and dilation, assume that for the
above mentioned coordinate system $i:D_M(u,\xi) \mapsto \reals^2$, we have
$i(u)=0$, $D(0,\rho) \subset i(D_M(u,\xi/2))$ for some
$\rho=\rho(u,\de)$ with  $0<\rho<\xi$
and if $x,x' \in D_M(u,\xi)$, then
\beq{amir-100}
(1-\de) |i(x)-i(x')| \leq d(x,x') \leq (1+\de) |i(x)-i(x')| \,.
\end{equation}

For any open $G \subseteq D_M(u,\xi)$, let
$\tau_G=\inf\{t \geq 0 : X_t \notin G\}$. It follows that for any $z\in
D_M(u,\xi)$
we can find a Brownian motion $B_t$ starting at $i(z)$ such that
$\{i(X_t),\,t\leq
\tau_G\}=\{B_{T_t},\,t\leq
\tau_G\}$ where $T_t=\int_0^t a(X_s)\,ds$,
see
\cite[Section V.1]{Revuz-Yor}.
Thus,
$T_{\tau_G}=\wt{\tau}_{i(G)}$, where  for any set $D\subseteq \reals^2$ we
write
$\wt{\tau}_{D}=\inf\{ t\geq 0 : B_t \notin D\}$. Consequently
\begin{equation}
\(\inf_{v\in G} a(v)\) {\tau_G}\leq  \int_0^{\tau_G}
a(X_s)\,ds= \wt{\tau}_{i(G)}\label{81.j}
\end{equation}

The upper bound in \req{mp.10} is obtained
by adapting the proof
provided in Section \ref{sec-hit}. To this end, fixing $1/2>\delta>0$,
extract a finite open sub-cover $\cup_{j} D_M(u_j,\xi_j/4)$ of the compact
manifold $M$ out of $\cup_{u\in M} D_M(u,\xi(u,\delta)/4)$.
Since $\uu{a} = \min_j \inf_{z \in D_M(u_j,\xi_j)} \; a_{u_j}(z)>0$,
we have by (\ref{amir-100}), (\ref{81.j}) and \req{amir-may1} that for any
$R\leq
\min_j \xi_j/4$
$$
\| \tau_R \|_R :=
\sup_{x \in M} \sup_{z \in D_M(x,R)} \E^z(\tau_{D_M(x,R)})
\leq \frac{R^2}{2 \uu{a} (1-\de)^2} \to_{R \to 0} \, 0.
$$
With its proof otherwise
unchanged,
Lemma \ref{lem-ld}
applies for $M$.
Moreover, fixing $j$, we have that for any
$x \in D_M(u_j,\xi_j/4)$ and $0<\ep<R<\xi_j/4$,
\beaa
i^{-1} ( D(i(x),(1-\de) \ep) ) & \subseteq & D_M(x,\ep)\, , \\
i^{-1} ( D(i(x),(1-\de) R) ) & \subseteq & D_M(x,R)\, , \\
i^{-1} ( D(i(x),(1-\de)^{-1} R/e) ) & \supseteq & D_M(x,R/e) \,.
\eeaa
Consequently, the left hand side of
\req{c2.12}
is bounded above by
the probability that $W_t$ does not hit
$D(i(x),(1-\de) \ep)$ during $n_\ep$ excursions, each starting at
$\partial D(i(x),(1-\de)^{-1} R/e)$ and ending at
$\partial D(i(x),(1-\de) R)$. This results with
\req{c2.12}
and hence
Lemma \ref{lem-hitprob}
holding, albeit with $1-\de'=(1-\de)(1+2\log(1-\de))$
instead of $(1-\de)$. Since $M$ is a smooth, compact, two-dimensional
manifold, there are at most $O(\ep^{-2})$ points $x_j \in M$ such that
$\inf_{\ell \neq j} d(x_\ell,x_j) \geq \ep$. The upper bound in \req{mp.10}
thus follows by the same argument that concludes
Section \ref{sec-upperbound}.

The complementary lower bound is next obtained
by adapting the proof
provided in Section \ref{sec-KRprth2}.
To this end, fixing $1/2>\de>0$,
let $\xi=\xi(\de)>0$ and $\rho=\rho(\de)>0$
be such that $D(0,\rho) \subset i(D_M(x_0,\xi/2))$ and
\req{amir-100} holds for the isothermal coordinate
system $i$ on $D_M(x_0,\xi)$, with
$i(x_0)=0$. It follows that
$$
\bigcup_{x \in S} D(x,\ep_1) \subset D(0,\rho) \subseteq i(D_M(x_0,\xi/2))\, ,
$$
provided $\ep_1<\rho/5$. Choosing $0<\ep_1<\rho/5$ small enough
so that $\ep_1 < R_1(\de)$ of Lemma \ref{lem-ld},
we say that $x \in S$ is $n$-successful if \req{pmperf} applies.
The probability $\oo{p}_k$ that
a planar Brownian path $B_t$ starting at
any $z \in \partial D(x,\ep_{n,k})$
hits $\partial D(x,\ep_{n,k-1})$ prior to
$\partial D(x,\ep_{n,k+1})$, is independent of
$z$ and
this is true
even after an arbitrary
random, path dependent, time change.
With $x_0 \notin i^{-1}(D(x,\ep_1))$,
and
$i^{-1}(\partial D(x,\ep_{n,k}))
\subset i^{-1}(D(0,\rho))$
for all $k=1,\ldots,n$,
we see that the identity \req{m1.1a}
holds, resulting with the conclusion of Lemma \ref{moment1}.
For $y \in S$, let $\GG_l^y$ be the $\sigma$-algebra generated by the
excursions $e^{(1)}, \ldots,e^{(N_{n,l}^y)}$, $e^{(N_{n,l}^y+1)}$
as defined in Section \ref{1momest}. Note that
Lemma \ref{lemprobends} applies
to the law of a planar Brownian excursion $B_\cdot$ starting at
$z \in \partial D(y,\ep_{n,l-1})$, conditioned to first exit
$D(y,\ep_{n,l})$ at $v$,
even after an arbitrary
random, path dependent, time change (indeed, both sides of
\req{m1.5eu} are clearly independent of such time change).
Moreover, the upper bound in \req{m1.5eu} is independent of the
initial point $z \in \partial D(y,\ep_{n,l-1})$.
In case $N^y_{n,l}=m_l>0$, since
$x_0 \notin i^{-1}(D(y,\ep_1))$ we have the representation
\req{m1.4e}, where conditioned upon $\GG^y_l$, the random vectors
$\{ Z^{(j)}_{k}, k \in J_l\}$ are independent for $j=1,2,\ldots,m_l$.
Recall the above mentioned identity between
the
`isomorphic image' of the path of $X_t$ till first exiting
$i^{-1}(D(y,\ep_{n,l}))$ and
the law of a time-changed planar Brownian path till its first
exit of $D(y,\ep_{n,l})$.
This identity, implies that
each random vector $\{ Z^{(j)}_{k}, k \in J_l\}$ has the conditional law of
$\{ Z_{k}, k \in J_l\}$ of Lemma \ref{lemprobends}
for some random $z_j \in \partial D(y,\ep_{n,l-1})$ and
$v_j \in \partial D(y,\ep_{n,l})$, both measurable on $\GG^y_l$.
With \req{m1.5eu} in force, we thus establish that the conclusion
\req{m1.2e} of Lemma \ref{recursionends} applies here
and can follow the proof of Lemma \ref{moment2} to arrive at its conclusion.
Thus establishing all estimates of Sections \ref{1momest}
and \ref{2momest}, we have that Lemma \ref{pmomlb} holds and
consequently the bound of \req{p-new2} applies.
It follows from \req{amir-100} that
\beaa
i^{-1}(\partial D(x,\ep_{n,n-1}))
&\subset& D_M(i^{-1}(x),(1+\de) \ep_{n,n-1})
\,,\\
 i^{-1}(\partial D(x,\ep_{n,n}))
&\bigcap& D_M(i^{-1}(x),(1-\de) \ep_{n,n}) = \emptyset\,,
\\
D_M(i^{-1}(x),(1-\de) \ep_{n,1})&\subseteq& i^{-1}(D(x,\ep_{n,1})) \,.
\eeaa
Consequently, if $x$ is $n$-successful, it follows that
$$
\TT(i^{-1}(x),(1-\de)\ep_{n,1}) \geq \sum_{j=0}^N \tau^{(j)} \,,
$$
where $N=n_n=3an^2\log n$ and $\tau^{(j)}$ correspond now to
excursions between the sets
$\partial D_M(i^{-1}(x),(1-\de)\ep_{n,n})$ and $\partial
D_M(i^{-1}(x),(1+\de)\ep_{n,n-1})$. The statement and proof of
Lemma \ref{lem1} then applies, except that we now use
$\TT (i^{-1}(x),(1-\de)\ep_{n,1})$ in \req{goalb}.
The lower bound in \req{mp.10} follows by the same argument as
in Section \ref{sec-KRprth2}, now with
$\CC_{(1-\de)\ep_{n(j),1}}$ in \req{p-new3}.
\qed

\section{Complements and unsolved problems}
\begin{enumerate}
\item 
We have the following direct corollary of Theorem \ref{theo-1p}.
\begin{corollary}\label{yuv-cor}
For $0< \ga < 1$ let $\TT_n(\ga)$ denote the time it takes 
until the largest disk unvisited by
the simple random walk in $\Z_n^2$ has radius $n^{\ga}$. 
Then,
$$
\lim_{n\to \ff}{\TT_n(\ga) \over (n\log n)^2}=
\frac{4(1-\ga)^2}{\pi}
\hspace{.1in}\mbox{in probability.}
$$
Equivalently, for $0<\al<1$ the logarithm to base $n$ of
the radius of the largest unvisited disk at time $\al \TT_n$
converges in probability to $1-\sqrt{\al}$.
\end{corollary}
\noindent{\bf Proof of Corollary \ref{yuv-cor}:}
The lower bound on $\TT_n(\gamma)$ is derived in Section \ref{sec-etc}.
By Theorem \ref{theo-1p} also 
$\PPP(\CC_{\ep_n} < \frac{2 (1-\ga+\de)^2}{\pi} (\log n)^2 ) \geq
1-\delta$ for $\ep_n = \frac{1}{3} n^{\gamma-1}$ 
and all $n$ large enough.
Similarly to Section \ref{sec-etc}, this yields
the upper bound on $\TT_n(\gamma)$ by strong approximation
(and tail estimates for  
the supremum of $|W_t-W_{k/2n^2}|$ over $t \in [k/2n^2,(k+1)/2n^2]$ and
$k \leq 4n^2 (\log n)^2$).
\qed
\item Given a planar lattice $\LL$, let
$\LL_\rho=\LL \cap D(0,\rho)$, a finite connected graph of $N_\rho$
vertices. Denote by $\TT_\rho$ the covering time for a
simple random walk on $\LL_\rho$. The approach of Section \ref{sec-etc}
can be adapted so as to show that
$$
\lim_{\rho \to \infty} \frac{\TT_\rho}{N_\rho (\log N_\rho)^2}
= C_\LL := \frac{A}{2\pi (\det \Gamma)^{1/2}} \qquad
\mbox { in probability }  \;,
$$
where
$$
A = \lim_{\rho \to \infty} \( \frac{\pi \rho^2}{N_\rho} \)
$$
is the area of a fundamental cell of $\LL$ and
$$
\Gamma = \lim_{n \to \infty} \frac{1}{n} \E ( S_n S_n' ) \;,
$$
is the two dimensional stationary covariance matrix associated with the
simple random walk on $\LL$ (note that $C_\LL$ is invariant
under affine transformations of $\reals^2$ and as such is an intrinsic
property of $\LL$). Of particular interest are the triangular (degree
$d=3$) and the honey-comb (degree $d=6$) lattices for which it is easy
to check that $\Gamma = \frac{1}{2} I$ and
$A=\frac{d}{4} \tan(\frac{\pi}{d})$.

\item Jonasson and Schramm show in \cite{Jonasson} the existence of
universal constants $C_d>0$ such that for any planar graphs $G_N$
of $N$ vertices and maximal degree $d_{\max} (G_N) \leq d$, one has
$$
\liminf_{N \to \infty} \frac{\TT(G_N)}{N (\log N)^2} \geq C_d \;,
$$
where $\TT(G_N)$ is the covering time for the simple random walk on
$G_N$. We believe that $C_d= \frac{d}{4\pi} \tan(\frac{\pi}{d})$ for
$d=3,4$ and $d=6$, corresponding to $G_N$ taken from the
triangular, square and honey-comb lattices, of degree $d=3,4$ and $6$,
respectively.

\item Recall that  $\TT_n$ denotes the (random)
cover time for simple random walk in  $\Z_n^2$.
A natural question, suggested to us by David Aldous, is to find
  a limit law for an appropriately normalized version of $\TT_n$.
The analogies with branching random walk lead us to supect that
perhaps the random variable $\TT_n^{1/2}/n$, minus its
median, will have a nondegenerate limit law.

\end{enumerate}

\vspace{2mm}
\noindent{\bf Acknowledgements}
We are grateful to David Aldous for
  suggesting the relevance of our results on
``thick points'' for random walks, to conjectures involving cover times.
We also thank  Isaac Chavel, Leon Karp, Mark Pinsky
and Rick Schoen for helpful discussions concerning
Brownian motion on manifolds.

\bigskip
\noindent
\begin{tabular}{lll}  & Amir Dembo& Yuval Peres\\
    & Departments of
Mathematics &  Departments of Mathematics \\
& and of Statistics & and of Statistics \\
   &Stanford University &   UC Berkeley\\
   &Stanford, CA 94305 &  Berkeley, CA 94720  \\
&amir@math.stanford.edu&peres@stat.berkeley.edu\\
& &\\
& & \\
& & \\
   & Jay Rosen & Ofer Zeitouni\\
    & Department of
Mathematics& Departments of EE and of Mathematics \\
   &College of Staten Island, CUNY& Technion, Haifa 32000, Israel\\
   &Staten Island, NY 10314& and Department of Mathematics, U. of Minnesota \\
&jrosen3@earthlink.net & Minneapolis, MN 55455 \\
&&zeitouni@math.umn.edu
\end{tabular}
\end{document}